\documentclass[leqno,11pt]{article}
\usepackage{latexsym}
\usepackage{amssymb}
\usepackage{amsfonts}
\usepackage{amsmath}
\usepackage{epsfig}
\usepackage{color}
\usepackage{pifont}
\usepackage{comment}
\usepackage{todonotes}
\usepackage{mathrsfs}
\usepackage{ulem}
\usepackage[unicode,breaklinks=true,colorlinks=true]{hyperref}

\allowdisplaybreaks
\makeatletter
\long\def\unmarkedfootnote#1{{\long\def\@makefntext##1{##1}\footnotetext{#1}}}
\makeatother

\newtheorem{definition}{Definition}[section]

\newtheorem{lemma}[definition]{Lemma}

\newtheorem{theorem}[definition]{Theorem}

\newtheorem{corollary}[definition]{Corollary}

\newtheorem{remark}[definition]{Remark}

\newtheorem{example}[definition]{Example}

\newcommand{\Le}{{\mathscr L}}

\def\o{\Omega}
\def\h{\mathcal H^n}
\def\H{\mathscr{H}}
\def\A{\mathcal A}

\def\vp{\varphi}

\def\m2{\frac{\h(\Omega)}2}
\def\tm2{\tfrac{\h(\Omega)}2}
\def\M2{\frac{|\Omega |}{2}}
\def\u+{u_+^*}

\def\-p{\overline{p}}

\def\w0{{W_0^{1,p}(\Omega)}}

\def\R{\mathbb R}
\def\N{\mathbb N}
\def\M{\mathbb M}

\def\ep{\varepsilon}

\newcommand{\supp}{\mathop{\rm supp}}

\def\rn{{{\R}^n}}

\newcommand{\Rn}{\mathbb{R}^{n}}
\newcommand{\Rm}{\mathbb{R}^{m}}
\newcommand{\LL}{\mathrm{L}}
\newcommand{\VV}{\mathrm{V}}
\newcommand{\WW}{\mathrm{W}}
\newcommand{\CC}{\mathrm{C}}
\newcommand{\dd}{\mathrm{d}}

\newcommand{\restrict}{\begin{picture}(10,8)\put(2,0){\line(0,1){7}}\put(1.8,0){\line(1,0){7}}\end{picture}}

\newcommand{\medint}{-\kern  -,395cm\int}
\newcommand{\medintinrigo}{-\kern  -,315cm\int}
\newcommand{\medelle}{-\kern  -,235cm L}
\newcommand{\medellenrigo}{-\kern  -,180cm L}
\newcommand{\qed}{\thinspace\null\nobreak\hfill
\hbox{\vbox{\kern-.2pt\hrule height.2pt
depth.2pt\kern-.2pt\kern-.2pt \hbox to1.8mm {\kern-.2pt\vrule
width.4pt \kern-.2pt\raise1.8mm\vbox to.2pt{} \lower0pt\vtop
to.2pt{}\hfil\kern-.2pt \vrule
width.4pt\kern-.2pt}\kern-.2pt\kern-.2pt \hrule height.2pt
depth.2pt \kern-.2pt}}\par\medbreak}

\title{Distortion of Hausdorff measures  under Orlicz-Sobolev maps}
\numberwithin{equation}{section}

\author{
Andrea Cianchi\\
\textit{Dipartimento di Matematica e Informatica \lq\lq U. Dini", Universit\`a di Firenze}\\
\textit{Viale Morgagni 67/A, 50134 Firenze, Italy} \\\textit{ \textrm{e-mail:} cianchi@unifi.it}
\bigskip
\\
Mikhail V.~Korobkov\\
\textit{School of Mathematical Sciences, Fudan University, Shanghai 200433, P.R.China}
\\
\textrm{and}
\\  
\textit{Sobolev Institute of Mathematics, pr-t Ac. Koptyug, 4, Novosibirsk, 630090, Russia}\\
\textit{ \textrm{e-mail:} korob@math.nsc.ru}
\bigskip
\\
Jan Kristensen\\
\textit{Mathematical Institute, University of Oxford, Andrew Wiles Building }
\\  
\textit{Oxford OX2 6GG, United Kingdom }
\\
\textit{ \textrm{e-mail:}  kristens@maths.ox.ac.uk }
}

\date{}

\begin{document}
\maketitle

 \begin{abstract} A comprehensive theory of the effect of Orlicz-Sobolev maps, between Euclidean spaces, on subsets with  zero or finite
Hausdorff measure is offered. Arbitrary Orlicz-Sobolev spaces embedded into the space of continuous function and Hausdorff measures built upon general gauge functions  are included in our discussion.
An explicit formula for the distortion of the relevant gauge function  under the action of these maps  is exhibited in terms of the Young function defining the Orlicz-Sobolev space.  New phenomena and features, related to the flexibility in the  definition of the degree of integrability of  weak derivatives of maps and in the notion of measure of sets, are detected.  Classical results, dealing with standard Sobolev spaces and Hausdorff measures, are recovered, and  their optimality is   shown to hold in a refined stronger sense. Special instances available in the literature, concerning Young functions and gauge functions of  non-power type, are also reproduced and, when not sharp,  improved. 
\end{abstract}

\pagestyle{myheadings} \thispagestyle{plain}

\unmarkedfootnote{\par\noindent\textit{Mathematics Subject
Classifications:} 46E35, 26B35.\par\noindent\textit{Key words and phrases: Hausdorff measure, gauge function, Orlicz-Sobolev space}.}

\section{Introduction and main results}\label{sec1}
A continuous map $u \colon \Omega \to \Rn$ in an open subset $\o$ of $\Rn$, with $n\in \N$, is said to satisfy the Lusin $N$ condition if 
it maps Lebesgue null sets into Lebesgue null sets. This amounts to requiring that, for sets $E\subset \Omega$,
\begin{equation}\label{lusin}
\text{if \,\,  $\Le^n (E) =0$,\,\, then \,\, $ \Le^n (u(E)) =0$,}
\end{equation}
where $\Le^n$ denotes the outer Lebesgue measure in $\Rn$.

The Lusin $N$ property plays a fundamental role in various results from classical real analysis and differentiation theory.  Early applications primarily concerning  the case $n=1$ can be found  in
\cite{Saks}.
Its crucial importance in   dimension $n>1$ is witnessed by its role for the validity of various formulas in  geometric measure theory 
and geometric function theory -- see e.g. \cite[Th.~3, p.~364]{RR} and \cite{MSZ1}. In particular, the Lusin $N$ condition is
critical in the proof of central properties of quasiconformal maps and, more generally, of 
maps with bounded distortion and of quasiregular maps. The latter were introduced by Yu.G.~Reshetnyak in ~\cite{Re0}  and subsequently
extensively analyzed and applied by several authors. A comprehensive treatment of the pertaining
theory and further references  can be found in the monographs \cite{Resh1,Vai,Ric1,Iw1}. 

 More recently a natural extension of the Lusin $N$ condition, where the Lebesgue measure is replaced
by lower dimensional Hausdorff measures, was decisive in works extending the Morse-Sard theorem and related results to classes of Sobolev maps
which  need neither  be Lipschitz nor everywhere differentiable -- see, in particular, \cite{bkk2013,bkk2015,KK18}. 

 In this connection, note that a map
$u \colon \Omega \to \Rm$, with $m \in \N$, which  is locally Lipschitz continuous, or everywhere differentiable, 
is easily seen to possess the property that, for
sets $E\subset \Omega$, 
\begin{equation}\label{haus}
\text{if \,\, $\H^\alpha (E) =0$, \,\, then \,\, $\H^\alpha (u(E)) =0$,}
\end{equation}
where $\alpha \in (0, n]$ and $\H^\alpha$ denotes the $\alpha$-dimensional Hausdorff measure.  In particular,   the Lusin $N$ property \eqref{lusin} follows from \eqref{haus}, when $m=n$. 

Information on the Hausdorff measure of the image of a set, which still implies the Lusin $N$ condition \eqref{lusin}, remains available if the map $u$,
though neither locally Lipschitz continuous nor everywhere differentiable, has yet a suitable degree of regularity.  With this regard, recall that the
space of locally Lipschitz continuous maps from $\o$ into $\Rm$ agrees with the  space $\VV^{1,\infty}_{\rm loc}(\o, \Rm )$.
Here, $\VV^{1,\infty}_{\rm loc}(\o, \Rm )$ stands for the local counterpart of  the homogeneous Sobolev space $\VV^{1,\infty}(\o, \Rm )$, where $\VV^{1,p}(\o, \Rm )$ is defined, for
$p \in [1, \infty]$, as
 \begin{equation}\label{sobolev}
  \VV^{1,p}(\o, \Rm ) = \big\{ u \in \LL^{1}_{\mathrm{loc}}( \o , \Rm )  \,\text{: $u$ is weakly differentiable and
    $|\nabla u| \in \LL^p (\o)$}\big\}.
\end{equation}
On extending earlier results from \cite{GV73} about dimension distortion under quasiconformal maps, in \cite{Kau}  Kaufman established  a modified version of property
\eqref{haus} for  maps in super-critical Sobolev spaces $\VV^{1,p}_{\rm loc}(\o, \Rm )$, with $p>n$. 
Recall  that the latter assumption ensures that each map $u \in \VV^{1,p}_{\rm loc}(\o, \Rm )$ has a (locally H\"older) 
continuous representative.
The result of \cite[Theorem 1]{Kau} asserts that, with this choice of the representative,
if $\alpha \in (0, n)$, then for sets $E\subset \o$,
\begin{equation}\label{kauf}
\text{ if\,  $\H^\alpha (E) < \infty$,  \,\, then \, $ \H^{\frac{\alpha p}{\alpha + p -n}} (u(E)) =0$.}
\end{equation}
The shift of the exponent from $\alpha$ to $\frac{\alpha p}{\alpha + p -n}$ is usually referred to as distortion of 
Hausdorff measure under the map $u$.
 Let us   notice that, by contrast, such a property does  not hold for all maps in $\VV^{1,p}_{\rm loc}(\o, \Rm )$ if $1\leq p \leq n$, even
under the additional assumption of continuity. With this regard, see \cite{C42,Po71,MM95} and, specifically, \cite[Theorem 4.3]{HeKo14}, for a streamlined presentation of Cesari's
example of a continuous map $u \in \VV^{1,n}(\Rn , \Rm )$, with $n \geq 2$ and $m \in \N$, such that $u\bigl( [-1,1] \times \{0 \}^{n-1} \bigr) = [-1,1]^m$.
This shows that the Lusin $N$ condition  fails in $\VV^{1,n}_{\rm loc}(\Omega, \Rn )$, for $n  \geq 2$. 
Let us also observe that the distortion property is neither
guaranteed by the mere membership of $u$ in the space of locally H\"older continuous maps with exponent $1-\frac np$, into which
$\VV^{1,p}_{\rm loc}(\o, \Rm )$ is embedded under the assumption that $p>n$. 




Several contributions on related questions are available in the literature -- see e.g.  \cite{ACTUTV13, bal2013, bkk2015, FKR, hencl-h2015, KKM99, KK14, KMZ15, Koskela, LSUT10, Raj}. 
Specifically,  \cite[Theorem C]{KKM99}  ensures that the Lusin $N$ property
holds for the continuous representative of every function from the Sobolev type space $\VV^{1}\LL^{n,1}_{\rm loc} (\o, \Rn)$ of those weakly differentiable maps
whose gradient belongs to the Lorentz space $\LL^{n,1}_{\rm loc}(\o, \Rn)$. This borderline Sobolev space, which coincides with $\VV^{1,1}_{\rm loc} (\o, \Rn)$ when $n=1$,
is  in fact optimal in the class  of all Sobolev type spaces built upon a rearrangement-invariant space for this property to hold.  The same space is
even optimal for embeddings into the space of continuous functions \cite{CianchiPick}. On the other hand,  the space $\VV^{1}\LL^{n,1}_{\rm loc}(\o, \Rn)$ is
not embedded into any space of locally uniformly continuous functions (with a common modulus of continuity), and no analogue of property \eqref{kauf}
holds for maps in $\VV^{1}\LL^{n,1}_{\rm loc}(\o, \Rn)$, since    they can map sets of any positive Hausdorff dimension into  sets of Hausdorff dimension~$n$. This means that they can map sets having zero $\alpha$-dimensional  Hausdorff measure, for any given $\alpha \in (0, n]$, into sets of positive  $n$-dimensional  Hausdorff measure.

The present paper is aimed at  
a comprehensive description of the  distortion of Hausdorff measures 
in the broader class of Orlicz-Sobolev maps.
These are maps in the homogeneous Orlicz-Sobolev space 
\begin{equation}\label{space}
 \VV^{1,A}(\o, \Rm ) = \big\{ u \in \mathrm{L}^{1}_{\mathrm{loc}}( \o, \Rm )  \,\text{: $u$ is weakly differentiable and
    $|\nabla u| \in \mathrm{L}^A(\o)$}\big \},
\end{equation}
or in its local counterpart  $\VV^{1,A}_{\rm loc}(\o, \Rm)$,
where $A$ is a Young function and $\mathrm{L}^A(\rn)$ is the associated Orlicz space.

In view of the discussion above, it is natural to focus on Orlicz-Sobolev spaces $\VV^{1,A}_{\rm loc}(\o, \Rm )$ which are embedded into a space of locally
uniformly continuous functions. A necessary and sufficient condition on $A$ for embeddings of this kind to hold, as well as the optimal
modulus of continuity are well known  \cite{Cianchi_ASNS}, and are recalled in Section \ref{prelim}. As a matter of fact, if $n\geq 2$ the same condition is also necessary and sufficient
for the embedding of $\VV^{1,A}_{\rm loc}(\o, \Rm )$ into the  space of plainly continuous functions, and even into $\LL^\infty_{\rm loc} (\o, \Rm )$.   In what follows, when dealing with a map from a space $\VV^{1,A}_{\rm loc}(\o, \Rm )$ associated with a Young function $A$ satisfying this  condition, we shall always refer to its continuous representative.

A treatment of the relevant questions,  tailored for the realm of general Orlicz-Sobolev spaces,  finds its proper environment in
a wider range of Hausdorff measures $\H^{\varphi}$, built upon  gauge functions $\varphi \colon [0, \infty) \to [0, \infty)$ which need not be just powers. In fact, this  extended
class of measures already appears in Hausdorff's original paper \cite{hausdorff}. The corresponding Hausdorff content, denoted by $\H^{\varphi}_\infty$,
will also come into play in our discussion.

Little seems to be known on the effect of Orlicz-Sobolev maps on general Hausdorff measures of  sets. The available results in this regard are apparently  limited to
the very special case when $A(t)\sim t^n (\log t)^q$ near infinity for some $q>0$, and $\varphi (r) \sim  r^\alpha (\log \frac 1r)^\beta$ near
$0$ for some $\beta >0$. They can be found in \cite{Koskela,KMZ15,Raj}.  Here, and in what follows, the relation
\lq\lq $\sim$" between two quantities means that they are bounded by each other, up to positive multiplicative constants.

Our contribution approaches the question in full generality and offers a handy recipe which, from any Young function $A$ and any gauge function $\varphi$ as above,
produces another gauge function $\psi$ with the property that if $u \in \VV^{1,A}_{\rm loc}(\o, \Rm )$, then, for sets $E\subset \o$,
\begin{equation}\label{ours}
\text{if \,\, $\H^\varphi (E) =0$\,\, then \,\, $ \H^{\psi} (u(E)) =0$.}
\end{equation}
Besides capturing  wider classes of Sobolev type maps and Hausdorff measures, broadening the functional framework entails the surfacing of certain  traits of the theory, which, by contrast, remain hidden for customary Sobolev spaces and   Hausdorff measures defined in terms of plain power type functions.     Allowing for a more flexible choice of function spaces and measures also enables us to supply sharper information about the optimality of well known results.

For instance, when $\varphi (r)=r^n$  property \eqref{ours} is shown to hold with $\psi (r)=r^n$. Hence, we infer that the Lusin $N$ property is supported
by any Orlicz-Sobolev space $\VV^{1,A}_{\rm loc}(\Omega, \Rn )$ which is embedded into a space of locally uniformly continuous functions.
This conclusion reproduces \cite[Corollary 2]{Kovt2014}, and can also be reached  via  Lusin $N$ property of the space $\VV^{1}_{\rm loc}\LL^{n,1}(\Omega, \Rn)$ recalled above,
owing to the inclusion $\LL^A_{\rm loc}(\o) \subset  \LL^{n,1}_{\rm loc}(\o)$, which is known to hold whenever
$\VV^{1,A}_{\rm loc}(\o, \Rn )$ is embedded into $\LL^\infty _{\rm loc}(\o, \Rn )$.

Even more interestingly,   the  comprehensive environment  of the present paper provides us with the required precision  to disclose a seemingly new  stability phenomenon, which is 
 not visible when the Sobolev spaces and the Hausdorff measures at disposal are the standard ones. It shows that
the validity of property \eqref{ours} with $\H^{\psi} \sim \H^\varphi$ is not peculiar to the Hausdorff measure $\H^n$, and hence, when $m=n$, to the
Lebesgue measure. 
Indeed, classes of Young functions $A$ and of gauge functions $\varphi$ are identified for which such an equivalence holds.
Roughly speaking, this happens when the function $A(t)$ has an essentially faster growth than the limiting power $t^n$ near
infinity -- e.g. $A(t)=t^p$ with $p>n$ -- and the gauge function $\varphi$  decays at zero  essentially more slowly than any power. One striking aspect
of this result is that it parallels the one holding for the Lebesgue measure, but at the opposite endpoint, in a sense, of the scale of Hausdorff
measures.  Conversely, for gauge functions $\varphi$ of  intermediate strength, or for Young functions $A$ with almost critical $n$-th power type
growth, a gap always occurs between the measures $\H^\varphi$ and $\H^{\psi}$.

Property \eqref{ours} will follow as a corollary of  Theorem \ref{main}. The latter amounts to a robust form of this property,
and tells us that, given any function $u\in \VV^{1,A}(\o, \Rm )$,  for every $\ep>0$ there exists $\delta>0$ such that,  for sets $E\subset \o$,
\begin{equation}\label{oursstable}
\text{if \,\, $\H^\varphi (E)<\delta$\,\, then \,\, $ \H^{\psi} (u(E)) <\ep$.}
\end{equation}
An analogue of property \eqref{oursstable} is also established with  Hausdorff measures replaced by Hausdorff contents.

These results are refined in Theorems  \ref{thm:measdist} and \ref{thm:contdist}.
The former provides us with a criterion on  the Young function $A$ and the gauge function $\varphi$ for an augmented version of property \eqref{ours} to hold.
This version  tells us that the vanishing assumption on $\H^\varphi (E)$ can just be replaced by a finiteness hypothesis, with the same conclusion. Namely,
 if $u \in \VV^{1,A}_{\rm loc}(\Omega, \Rm )$, then, for sets $E\subset \o$,
\begin{equation}\label{oursaug}
\text{if \,\, $\H^\varphi (E)< \infty$\,\, then \,\, $ \H^{\psi} (u(E)) =0$.}
\end{equation}
The examples alluded to above, when $\H^\varphi \sim \H^{\psi}$, demonstrate that property \eqref{oursaug} fails
if $A$ and $\varphi$ do not meet the requirements of the relevant criterion. However, in this case 
a quantitative form of  property  \eqref{ours} is yet offered, which tells us that 
\begin{equation}\label{distintro}
\H^{\psi}\bigl( u(E) \bigr) \leq  c \,\H^{\vp}(E)
\end{equation}
for   maps $u \in \VV^{1,A}(\Omega, \Rm )$, for all sets $E \subset \Omega$ and
for some constant $c$ depending, in an explicit way, on $n$, $\psi$ and $\| \nabla u \|_{\LL^{A}( \Omega )}$. 

Finally, Theorem \ref{thm:contdist} contains a variant of inequality \eqref{distintro} involving Hausdorff contents instead of measures.
The inequality in question holds for  maps in $\VV^{1,A}(\Rn, \Rm )$. Moreover, the dependence on $\H^{\vp}_\infty(E)$ on the right-hand side is nonlinear.
A slightly more general  version of this inequality can be established, for maps defined in  arbitrary open sets $\Omega \subset \Rn$, in the case
when $\H^{\vp}_\infty(E)$ agrees with the outer Lebesgue measure. This result is enucleated in Theorem \ref{nlusin}.

Precise statements of the    results outlined above are given in Section \ref{S:main}, where applications  to special choices of the functions $A$
and $\varphi$ are also proposed. In particular, they recover Kaufman's theorem \eqref{kauf},
improve a theorem of \cite{Raj},  and specify the results highlighted   above. Section \ref{tech} is devoted to  preliminary technical lemmas. The proofs  of the main results are then accomplished in Section \ref{proofs}.
In the final Section \ref{S:cex}, the sharpness of some of the examples
worked out is  discussed. In particular, as hinted above, thanks to the use of general Hausdorff measures, an enhanced optimality of the
results already available in the literature can be established.

\section{Hausdorff measures and Orlicz-Sobolev spaces}\label{prelim}

In this section we collect   a few definitions and properties about the Hausdorff measures and the  underlying gauge functions considered
throughout, and the  Young functions for the admissible Orlicz-Sobolev spaces in our main result.

\subsection{Gauge functions and Hausdorff measures}\label{sub1}

A Hausdorff measure can be built upon any gauge function, an expression whereby we understand any function $\varphi \colon [0, \infty) \to [0, \infty)$ which is increasing,  continuous
and  vanishes only at $0$.  Set, for $\sigma \in (0, \infty]$, 
\begin{equation}\label{jan40}
\H_\sigma^\varphi (E) = \inf\left\{\sum_{i=1}^\infty\varphi\bigl(d(E_i)\bigr)\ :\ E_i\subset \Rn,\,\,d(E_i)\le \sigma,\ \ E \subset\bigcup\limits_{i=1}^\infty E_i \right\}
\end{equation}
for any set $E\subset \Rn$. Here, $d(\cdot)$ stands for diameter. 
Then 
the relevant measure  $\H^{\varphi}$   is  defined as 
\begin{equation}\label{Hmeas}
\H^{\vp} (E)= \lim _{\sigma \to 0^+} \H_\sigma^\varphi (E)
\end{equation}
The set function $\H_\infty^\varphi$  also enters our results and is called Hausdorff content.
\\ For sets $E\subset \Rn$, one has that
\begin{equation}\label{jan60}
\H_\infty^\varphi (E) = 0 \quad \text{if and only if} \quad \H^\varphi (E)=0,
\end{equation}
see  \cite[Proposition 5.1.5]{AH}.
\\
In order to rule out trivial cases, we shall assume that
\begin{equation}\label{feb3}
\liminf_{r\to 0^+}\frac{\varphi (r)}{r^n} >0.
\end{equation}
Actually, if condition \eqref{feb3}  fails, then $\H^\varphi (E) =0$ for every set $E\subset \Rn$. This follows from \cite[Proposition 5.1.8]{AH}.
In fact, a minor variant of this proposition also ensures that, under \eqref{feb3}, one can additionally assume, without loss of generality, that the function 
\begin{equation}\label{feb8}
\text{$\frac{\varphi (r)}{r^n}$ is non-increasing,}
\end{equation}
see Lemma \ref{adamshedberg518} below.
\\
We shall therefore suppose throughout,  without explicit mentioning, that all gauge functions $\varphi$  fulfill condition   \eqref{feb8}.
\\ Clearly,   choosing $\varphi (t)=t^\alpha$, with $\alpha \in (0,n]$,  in the definitions of  $\H^\varphi$ and $\H_\infty^\varphi$ reproduces (up to a multiplicative constant) the
customary Hausdorff measure $\H^{\alpha}$ and Hausdorff content  $\H^{\alpha}_\infty$.
\\ One has that $\H^{\varphi}\neq \mathcal L^n$ is and only if
\begin{equation}\label{Feb4}
\lim_{r\to 0^+}\frac{\varphi (r)}{r^n}=\infty.
\end{equation}
Observe that, thanks to property \eqref{feb8},
\begin{equation}\label{delta2}
\varphi (k r) \leq \max\{1, k^n\} \,\varphi (r) \qquad \text{for $k>0$ and $r \geq 0$.}
\end{equation}
An alternate Hausdorff type measure  $\Lambda^{\varphi}$, defined in terms of coverings by dyadic closed cubes, will be needed in our proofs. These measures are sometimes called net measures, see
\cite[Chap.~5]{Fa85}. On setting
\begin{equation}\label{jan41}
\Lambda _\sigma ^\varphi (E) =   \inf \bigg\{\sum _{i}\varphi (d(Q_i)): E \subset \bigcup_i Q_i, \,\, Q_i \,\text{is a dyadic cube and }
  d(Q_{i}) \leq \sigma \bigg\}
\end{equation}
for any set $E\subset \Rn$, the relevant measure
 is defined as 
\begin{equation}\label{Lambda}
\Lambda^{\vp} (E)= \lim _{\sigma \to 0^+}  \Lambda _\sigma ^\varphi (E).
\end{equation}
Thanks to property 
\eqref{delta2},   
\begin{equation}\label{apr1}
\H^{\varphi}_\sigma(E) \leq \Lambda^{\varphi}_{\sigma}(E) \leq c_{n}\H^{\vp}_{\frac \sigma 2}(E)
\end{equation}
for $\sigma \in (0, \infty]$ and for every $E \subset \rn$, and hence
\begin{equation}\label{feb25}
\H^{\varphi}(E) \leq \Lambda^{\varphi}(E) \leq c_{n}\H^{\vp}(E)
\end{equation}
for every set $E \subset \rn$, where the constant $c_{n}$ only depends on the
dimension $n$ (see, e.g., \cite[Lemma~2.2]{bkk2015}\,). 
Moreover,
by property \eqref{delta2} again, given any $k>0$, on defining  the gauge function $\varphi_k$ as 
\begin{equation}\label{phik}
\varphi_k (t) = \varphi (kt) \quad \text{ for $t\geq 0$,}
\end{equation}
 one has that
\begin{equation}\label{feb26}
\Lambda ^{\varphi_k }(E) \approx  \Lambda^{\varphi} (E)
\end{equation}
for   every set $E \subset \rn$, up to multiplicative constants depending only on $k$. 
\\
Since $\H^{\varphi}$ is a regular outer measure, it enjoys the ascending sets property, which tells us that, if $\{E_j\}$ is any sequence of sets in $\rn$ such that
$E_j\subset E_{j+1}$ for $j\in \N$, then 
\begin{equation}\label{subad}
\H^{\varphi}\Big(\bigcup_{j}E_j\Big) = \lim_{j \to \infty} \H^{\varphi}(E_j),
\end{equation}
see  \cite{Ro70/98}.
A deeper result from \cite{Da70} (see also \cite{Ro70/98}) ensures that the same conclusion holds for the Hausdorff content.  Namely, if $\{E_j\}$  is a sequence as above, then
\begin{equation}\label{subcont}
\H^{\varphi}_\infty\Big(\bigcup_{j}E_j\Big) = \lim_{j \to \infty}  \H^{\varphi}_\infty(E_j).
\end{equation}
Notice that we shall only make use of \eqref{subcont} in the case when the gauge function $\vp$ is a power. In fact, we could, at the expense of an additional multiplicative constant, use the much more elementary result \cite[Lemma 5.3]{Fa85}. Instead of \eqref{subcont}, the latter  yields  the bound
$$
\H^{\varphi}_\infty\Big(\bigcup_{j}E_j\Big) \leq  c\lim_{j \to \infty}  \H^{\varphi}_\infty(E_j)
$$
for some constant that only depends on $\vp$ and $n$.
\begin{lemma}\label{adamshedberg518}
Let $\varphi \colon [0, \infty ) \to [0, \infty )$ be a continuous increasing function such that $\varphi (0)=0$ and satisfying condition \eqref{feb3}.
Then, there exist a continuous  increasing function $\varphi ^{\circ}\colon [0,\infty ) \to [0,\infty )$ and a positive constant  $c=c(n)$  such that $\varphi ^{\circ}(0)=0$, 
\begin{equation}\label{july200}
\text{$\frac{\varphi ^{\circ}(r)}{r^n}$ is non-increasing,}
\end{equation}
and
\begin{equation}\label{equivhaus}
c\,\H^{\varphi }_{\delta}(E) \leq \H^{\varphi ^{\circ}}_{\delta}(E) \leq \H^{ \varphi }_{\delta}(E)
\end{equation}
for every $\delta \in (0,\infty ]$ and   every set $E \subset \Rn$.
\end{lemma}

\noindent
\textit{Proof.}
Define the function  $\varphi ^{\circ}\colon [0,\infty ) \to [0,\infty )$ as
$$
\varphi^{\circ}(r) = \left\{
\begin{array}{ll}
\displaystyle{r^{n}\inf_{0<t \leq r} \frac{\varphi(t)}{t^n}} & \mbox{ if } r > 0\\
  0 & \mbox{ if } r = 0.
\end{array}
\right.
$$
One can verify that the function  $\varphi^\circ$ is increasing and fulfills  property \eqref{july200}. These monotonicity properties imply that $\varphi^\circ$ is continuous in $(0, \infty)$, and, since $\varphi$ is continuous in $[0, \infty)$ and vanishes at $0$,  the function $\varphi^\circ$ is continuous in $[0, \infty)$ as well.
\\ 
Inasmuch as  $\varphi^{\circ} \leq \varphi$,  the second inequality  in equation \eqref{equivhaus} is trivial. In order to prove  the first one,  fix $\delta \in (0, \infty ]$ and $E \subset \Rn$. We may assume that $\H^{\varphi^{\circ}}_{\delta}(E) < \infty$, 
otherwise there is nothing to prove. Given any $\varepsilon > 0$,  pick a covering $ \bigcup_{i=1}^{\infty}E_i \supset E$, with $0< d(E_{i}) < \delta$ and such that
\begin{equation}\label{july204}
\sum_{i=1}^{\infty}\varphi^\circ\bigl( d(E_{i}) \bigr) < \H^{\varphi^\circ}_{\delta}(E)+\varepsilon .
\end{equation}
For each $i$ we choose $d_{i} \in (0,d(E_{i})]$ with the property that $\varphi^{\circ}(d(E_{i}))(1+\varepsilon ) > \bigl( d(E_{i})/d_i \bigr)^{n}\varphi(d_{i})$. Now, observe that  there exists a dimensional constant $c$ such that, given $0<d\leq D$, any subset of $\Rn$ of diameter
$D$  can be covered by a family of sets with cardinality less that  $c\bigl( D/d \bigr)^n$, and   diameter at most $d$. This fact can be verified by considering 
a suitable translation and dilation of the grid defined by $\mathbb{Z}^n$. 
Hence, for each $i$ the set $E_i$ can be decomposed as $E_i = \bigcup_{j} E_{i,j}$ where the cardinality of the family of the sets 
 $E_{i,j}$ does not exceed $c\bigl( d(E_{i})/d_{i}\bigr)^n$,  and $d(E_{i,j}) \leq d_i$. Consequently,
\begin{equation}\label{july203}
  (1+\varepsilon )\sum_{i=1}^{\infty} \varphi^{\circ}\bigl( d(E_{i}) \bigr) > \sum_{i=1}^{\infty} \left( \frac{d(E_{i})}{d_{i}} \right)^{n}\varphi(d_{i})  \geq \frac{1}{c}\sum_{i=1}^{\infty} \sum_{j} \varphi \bigl( d(E_{i,j}) \bigr) \geq \frac{1}{c}\H^{\varphi}_{\delta}(E).
\end{equation}
The conclusion follows from inequalities \eqref{july204} and \eqref{july203}. \qed

\smallskip

\subsection{Young functions and Orlicz-Sobolev spaces}\label{sub2}

A Young function $A \colon [0, \infty) \to [0, \infty]$ is a (non trivial)  convex function such that $A(0)=0$. 
As a consequence, the function 
\begin{equation}\label{monotone}
\text{$\frac{A(t)}t$ is non-decreasing in $(0,\infty)$.}
\end{equation}
\\
The Young conjugate $\widetilde{A}$ of $A$  is also a Young function and is defined by
\begin{equation}\label{2.8}
\widetilde{A}(t) = \sup \bigl\{\tau t-A(\tau ):\,\tau \geq 0 \bigr\}  \quad \text{for $t\geq 0$.}
\end{equation}
One has that $\widetilde{\widetilde A} = A$.
\\ A  Young function $A$ is said to dominate another Young function $A_0$ globally 
   if there exists a positive
 constant $c$  such that
\begin{equation}\label{B.5bis}
A_0(t)\leq A(c t) \quad \text{for $ t\geq 0.$}
\end{equation}
The function $A$ is said to dominate $A_0$ near infinity  if there
exists $t_0> 0$ such that \eqref{B.5bis} holds for $t \geq t_0$.
The functions $A$ and $A_0$ are called equivalent globally [near infinity]    if they dominate each other globally [near infinity].

Let $\Omega$ be a measurable set in $\Rn$ and let $A$ be a Young function. The Orlicz space $\LL^A (\Omega)$    is the Banach
space of measurable functions $u\colon \Omega \to \R$ for which the Luxemburg norm
\begin{equation}\label{lux}
 \|u\|_{\LL^A(\Omega)}= \inf \left\{ \lambda >0 :  \int_{\Omega }A
\left( \frac{|u|}{\lambda} \right) dx \leq 1 \right\}\,
\end{equation}
is finite. 
\\ If a Young function $A$ dominates another Young function $A_0$ globally [near infinity], then 
\begin{equation}\label{jan10}
\LL^A (\Omega) \hookrightarrow \LL^{A_0} (\Omega)
\end{equation}
for every measurable set $\Omega$ [for every measurable set $\Omega$  with $\Le^n (\Omega)<\infty$]. Moreover, embedding \eqref{jan10} holds as an equality (up to equivalent norms), if equivalence holds between the Young function $A$ and $A_0$ instead of domination.
\\  When $\Omega$ is an open set, the corresponding homogeneous Orlicz-Sobolev space  $\VV^{1,A}(\Omega, \Rm )$ is defined by \eqref{space}.
The local Orlicz-Sobolev space  $\VV^{1,A}_{\rm loc}(\Omega, \Rm )$ is defined accordingly.

As explained in Section \ref{sec1}, our analysis applies to Orlicz-Sobolev spaces $\VV^{1,A}_{\rm loc}(\Omega, \Rm )$ supporting  an embedding of the form 
\begin{equation}\label{embeddC}
\VV^{1,A}_{\rm loc}(\Omega, \Rm ) \hookrightarrow \CC^\omega_{\rm loc} (\Omega, \Rm)
\end{equation}
 for some space of (locally) uniformly continuous functions
$\CC^\omega _{\rm loc}(\Rn, \Rm)$, with modulus of continuity  $\omega \colon [0, \infty) \to [0, \infty)$.
By  \cite[Theorems 1a and 3]{Cianchi_ASNS}, such an embedding holds
 if and only if
\begin{equation}\label{feb1}
\begin{cases}
\displaystyle \lim_{t\to \infty} \frac t{A(t)} = 0 & \quad \text{if $n=1$,}
\\ \\ \displaystyle
\int^\infty \bigg(\frac t{A(t)}\bigg)^{\frac 1{n-1}}\, dt < \infty & \quad \text{if $n\geq 2$.}
\end{cases}
\end{equation}
When $n\geq 2$ the condition in \eqref{feb1} is in fact also necessary and sufficient for the embedding of $\VV^{1,A}_{\rm loc}(\Omega, \Rm )$
into the  space of plainly continuous functions, and even into $\LL^\infty_{\rm loc} (\Omega, \Rm )$. The sufficiency of this condition was  earlier established  in ~\cite{Calderon}.
\\
Specifically, in connection with \eqref{embeddC},
 denote by $B$ the Young function given by
\begin{equation}\label{B}
B(t) = \begin{cases}
A(t) & \quad \text{if $n=1$}
\\ \\
 \Bigg( t^{n'}\displaystyle \int_t^\infty \frac{\widetilde A(s)}{s^{1+n'}}\, ds\Bigg)^{\widetilde {\quad}} & \quad \text{if $n \geq 2$,}
\end{cases}
\end{equation}
where the notation $(\dots)^{\widetilde {\quad}}$ stands for the Young conjugate of the function in parentheses.
Observe that, by \cite[Lemma 4.4]{Cianchi_Ibero}, the convergence of the integral in \eqref{B} is equivalent to that in \eqref{feb1}.
Then
there exists a constant $\kappa = \kappa (n,m)$ such that the Morrey type inequality
\begin{equation}\label{poincare}
\frac{d(Q)^n}\kappa B\bigg(\frac {d(u(Q))}{\kappa \,\lambda \, d(Q)}\bigg) \leq \int_Q A\bigg(\frac{|\nabla u|}\lambda\bigg)\, dx,
\end{equation}
holds for the continuous representative of any map in $u\in \VV^{1,A}_{\rm loc}(\Omega, \Rm )$. Here $Q$ denotes any cube compactly contained in 
$\Omega$  and $\lambda$ any positive number which makes the integral finite -- see \cite[Proof of Theorem 1.1]{albericocianchi}
(and also \cite[Theorem 4.1]{carozzacianchi}). Moreover,  if  $u\in \VV^{1,A}(\Omega, \Rm )$, inequality \eqref{poincare}   holds also for every open cube $Q \subset \Omega$.
\\
In inequality \eqref{poincare},  the continuous representative of $u$ is considered. 
As already mentioned above,  such a  representative will always  be considered also in what follows, without explicitly mentioning,
when dealing with maps in the space $ \VV^{1,A}_{\rm loc}(\Omega, \Rm )$ with a Young function $A$ fulfilling assumption \eqref{feb1}.


In order to avoid additional unnecessary technicalities, in our discussion of the distortion of Hausdorff measures  we shall  assume that the Young function $A$ defining the ambient Orlicz-Sobolev space is such that
\begin{equation}\label{0inf}
0<A(t)<\infty \quad \text{for }t>0,
\end{equation}
and
\begin{equation}\label{feb7}
\begin{cases}
\displaystyle \lim_{t\to 0^+} \frac t{A(t)} = \infty & \quad \text{if $n=1$,}
\\ \\ \displaystyle
\int_0 \bigg(\frac t{A(t)}\bigg)^{\frac 1{n-1}}\, dt = \infty & \quad \text{if $n\geq 2$.}
\end{cases}
\end{equation}
A few comments on  assumptions \eqref{0inf} and \eqref{feb7} are in order.
The assumptions that $A(t)>0$ for $t>0$ and that $A$ satisfies \eqref{0inf} are immaterial when dealing with the space  $\VV^{1,A}_{\rm loc}(\Omega, \Rm)$. The same is true for  the space  $\VV^{1,A}(\Omega, \Rm )$, provided that $\Le^{n} (\Omega)<\infty$. 
Indeed, the function $A$ can be replaced, if necessary,  by a Young function $A_0$ equivalent near infinity such that  $A_0(t)>0$ for $t>0$ and satisfying condition \eqref{feb7}.
Thanks to  equation \eqref{jan10}, after this replacement,   $\VV^{1,A}_{\rm loc}(\Omega, \Rm )=  \VV^{1, A_0}_{\rm loc}(\Omega, \Rm )$, and,
if $\Le^{n} (\Omega)<\infty$, then $\VV^{1,A}(\Omega, \Rm ) = \VV^{1, A_0}(\Omega, \Rm )$  (up to equivalent norms). By contrast, such a replacement does affect the space $\VV^{1,A}(\Omega, \Rm )$ when $\Le^{n} (\Omega)=\infty$. However, a close inspection of our proofs will reveal  versions of our main results stated for maps in the space $\VV^{1,A}(\Omega, \Rm )$ continue to hold  in the broader   space 
of those weakly differentiable maps $u : \Omega \to \Rm$ such that  $|\nabla u| \in \mathrm{L}^A(E)$ for every set $E\subset \Omega$ such that $\Le^{n} (E)<\infty$. The latter space is again invariant under the replacement of $A$ by a function $A_0$  as above.
\\ As for the requirement that $A(t)<\infty$ for $t >0$, it can actually been dropped   in our results about properties \eqref{ours} and \eqref{oursstable}, and in their counterparts about Hausdorff contents. This is due to the fact that, if $A(t)=\infty$ for large values of $t$, then $\VV^{1,A}(\Omega, \Rm ) \hookrightarrow \VV^{1,\infty}(\Omega, \Rm )$. Since the proofs of the properties in question only make use of the local behaviour of maps in $\VV^{1,A}(\Omega, \Rm )$, our results follow via the classical results for Lipschitz continuous maps recalled above, inasmuch as maps in $\VV^{1,A}(\Omega, \Rm )$ are locally Lipschitz continuous. On the other hand, the case of Lipschitz continuous functions demonstrates that the finiteness assumption on $A$ cannot be dispensed with when the stronger property \eqref{oursaug} is in question.

%
%
%

Let us finally observe that,  when $n \geq 2$,  the function $A$ dominates  globally the function $B$    given by \eqref{B} \cite[Lemma 3.2]{carozzacianchi}.
Moreover, the functions $A$ and $B$ are equivalent near infinity if and only if 
\begin{equation} \label{indexcond}
i_\infty(A) >n,
\end{equation}
where  $i_\infty(A)$ denotes the lower Matuszewska--Orlicz  index of $A$ at infinity.  For Young functions $A$ fulfilling condition \eqref{0inf}, this index is defined as 
\begin{equation}\label{index}
i_\infty(A)= 
\lim _{\lambda \to \infty} \frac{\log \Big(\liminf _{t\to
\infty }\tfrac{A(\lambda t)}{A(t)}\Big)}{\log \lambda }.
\end{equation} 
The functions $A$ and $B$ are equivalent globally if and only if
\begin{equation} \label{indexcondglob}
i(A) >n\,
\end{equation}
where  $i(A)$ denotes the global lower   Matuszewska–Orlicz  index of $A$, given by 
\begin{equation}\label{indexglob}
i(A)= 
\lim _{\lambda \to \infty} \frac{\log \Big(\inf _{t > 0}\tfrac{A(\lambda t)}{A(t)}\Big)}{\log \lambda }.
\end{equation} 
See \cite[Lemma 2.3]{Strom} and, for more details, \cite[Proposition 4.1]{CianchiMusil}.

\section{Main results and examples}\label{S:main}
The gauge function $\psi$ which renders implication \eqref{ours} true is associated with any Young function $A$ satisfying assumptions
\eqref{feb1}, \eqref{0inf}, \eqref{feb7}, and any   gauge function $\varphi$  as follows.
Consider the function $J \colon [0, \infty) \to (0, \infty)$ defined as
\begin{equation}\label{J}
J(s) =\begin{cases} s\, B^{-1}\Big(\displaystyle\frac {\varphi (s)}{s^n}\Big) & \quad  \text{for $s>0$}
\\ 0 & \quad  \text{for $s=0$,}
\end{cases}
\end{equation}
where $B$ is defined at \eqref{B}.
Then the function $\psi \colon [0, \infty) \to [0, \infty)$ is given by
\begin{equation}\label{psi}
\psi (r) =   \varphi(J^{-1}(r)) \qquad \text{ for $r>0$.}
%
\end{equation} 
The   function $\psi$ will be shown to be a well-defined gauge function in the next section. 

Let us   point out  that, in view of the definition of the Hausdorff measures $\H^\varphi$ and $\H^\psi$, only the asymptotic behavior of $B$,
and hence of $A$, near infinity, and the asymptotic behaviours of $\varphi$, $J$ and $\psi$ near zero are relevant.
In particular, if $A$ fulfills property \eqref{indexcond}, then, owing to inequality \eqref{delta2psi} below, replacing $B$ simply by $A$ in definitions \eqref{J}
and \eqref{psi} results in an equivalent  function $\psi$, and hence in an equivalent Hausdorff measure $\H^\psi$.
Analogously, under assumption \eqref{indexcondglob}, the Hausdorff content $\H_\infty^\psi$ can be equivalently defined just via $A$ instead of $B$.

\smallskip

We are now in a position to state our first main result, which concerns property \eqref{oursstable}.

\begin{theorem}\label{main}
Let $n$, $m \in \N$, let $\Omega$ be an open set in $\R^n$, and let $A$ be a Young function satisfying   conditions \eqref{feb1}, \eqref{0inf} and \eqref{feb7}.  Assume that $\varphi$ is a  gauge
function, and let $\psi$ be the gauge function defined by \eqref{psi}. 
\\ (i) Let $u\in \VV^{1,A}(\Omega, \Rm )$. For every $\ep >0$ there exists $\delta >0$ such that, for sets $E \subset \Omega$,
\begin{equation}\label{feb6}
\text{if $\H^{\varphi }(E) < \delta$, \quad then \quad $\H^{\psi}(u(E))<\ep$.}
\end{equation}
\\ (ii) Let $u\in \VV^{1,A}(\Rn, \Rm )$.  For every $\ep >0$ there exists $\delta >0$ such that, for sets $E \subset \Rn$,
\begin{equation}\label{Feb60}
\text{if $\H^{\varphi }_\infty(E) < \delta$, \quad then \quad $\H^{\psi}_\infty(u(E))<\ep$.}
\end{equation}
\end{theorem}

\begin{remark}\label{NforOrlicz}{\rm
Observe that, if $\varphi (r)=r^n$, then $\psi (r)=r^n$ for any Young function $A$ as in the statement of Theorem \ref{main}. Hence, one infers for $u \in \VV^{1,A}( \Omega , \Rn )$ that
 for every $\ep >0$ there exists $\delta >0$ such that, for sets $E \subset \Omega$,
\begin{equation}\label{dic50}
\text{if $\Le^n(E) < \delta$ \quad then \quad $\Le^n(u(E))<\ep$.}
\end{equation}
Thus the classical Lusin $N$-property   for maps in the Orlicz-Sobolev space $\VV^{1,A}(\Omega, \Rn )$, with $A$ fulfilling condition \eqref{feb1},
announced  in Section~\ref{sec1} is in particular recovered from Theorem~\ref{main}.  }
\end{remark}
 
Property \eqref{ours} is the content of the next corollary and follows from Theorem \ref{main}.

\begin{corollary}\label{cor}
Let $n$, $m \in \N$, let $\Omega$ be an open set in $\R^n$, and let $A$ be a Young function satisfying   conditions  \eqref{feb1}, \eqref{0inf} and \eqref{feb7}.  Assume that $\varphi$ is a  gauge
function, and let $\psi$ be the gauge function defined by \eqref{psi}. 
Let $u\in \VV^{1,A}_{\rm loc}(\Omega, \Rm )$. Then for sets $E \subset \Omega$,
\begin{equation}\label{cor1}
\text{if $\H^{\varphi }(E)=0$, \quad then \quad $\H^{\psi}(u(E))=0$.}
\end{equation}
\end{corollary}

Formulating the qualifications on the Young function $A$ and the gauge function $\varphi$ which entail the validity of the enhanced property  \eqref{oursaug} involves the  \lq\lq scaling functions" defined as follows.
\\
With any  function $h \colon (0, \infty ) \to (0, \infty )$,  we associate the 
function
 $\Theta_{h}^0 \colon (0,\infty ) \to [0,\infty ]$,  defined as
 \begin{equation}\label{defthetah}
\Theta_{h}^0 (r) =  \limsup_{t \to 0^{+}} \frac{h(rt)}{h(t)} \quad \text{for $r>0$,}
\end{equation}
the function
$\Theta_{h}^\infty \colon (0,\infty ) \to [0,\infty ]$, defined as 
 \begin{equation}\label{defthetaH}
\Theta_{h}^\infty (r) =  \limsup_{t \to \infty} \frac{h(rt)}{h(t)} \quad \text{for $r>0$,}
\end{equation}
and the function $\Theta_{h} \colon (0,\infty ) \to [0,\infty ]$, defined as
\begin{equation}\label{defthetahs}
\Theta_{h}(r) = \sup_{t > 0} \frac{h(rt)}{h(t)} \quad \text{for $r>0$.}
\end{equation}
In particular, if the function $h$ is non-decreasing, then these three functions are non-decreasing as well.
They  naturally arise in connection with scaling properties of 
Hausdorff measures and Hausdorff contents. 

 
\begin{theorem}\label{thm:measdist}
Let $n$, $m \in \N$ and let $\Omega$ be an open subset of $\Rn$. Let $A$ be a Young function satisfying
conditions \eqref{feb1}, \eqref{0inf}, \eqref{feb7}, and let $B$ be the function defined by \eqref{B}.   Assume that $\varphi$ is a   gauge
function fulfilling condition  \eqref{Feb4}, and 
let $\psi$ be the gauge function defined by \eqref{psi}. 
\\ (i) Assume that 
\begin{equation}\label{scalingassump}
\text{
$\lim_{r \to 0^{+}} \Theta_{\psi}^0(r)=0$ \quad and \quad $\lim_{r \to  0^{+}} \Theta_{B^{-1}}^\infty(r) =0$. }
%
\end{equation}
 Let $u\in \VV^{1,A}_{\rm loc}(\Omega, \Rm )$.  Then, for  sets  $E\subset \Omega$,
\begin{equation}\label{Feb1}
\text{if $\H^\vp(E)<\infty$, \quad then \quad  $\H^{\psi}\bigl( u(E) \bigr) =0$.}
\end{equation}
(ii) Assume that 
\begin{equation}\label{scalingassump1}
\text{either \quad 
$\lim_{r \to 0^{+}} \Theta_{\psi}^0(r)>0$ \quad or \quad $\lim_{r \to  0^{+}} \Theta_{B^{-1}}^\infty(r) >0$. }
%
\end{equation}
 Let $u \in\mathrm{V}^{1,A}(\Omega, \Rm )$.   Then 
\begin{equation}\label{measdistort}
\H^{\psi}\bigl( u(E) \bigr) \leq c \, \Theta_{\psi}^0\bigl( \kappa \| \nabla u \|_{\LL^{A}( \Omega )}^{+} \bigr) \H^{\vp}(E)
\end{equation}
for every set  $E\subset \Omega$, where $\kappa$ is the constant appearing in inequality \eqref{poincare},
$c= c_n \lim_{r \to 0^{+}}   \Theta_{\psi}^0\bigl( \Theta_{B^{-1}}^\infty(r)\bigr)  $, and $c_n$ is the constant appearing in equation \eqref{feb25}. On the right-hand side of \eqref{measdistort}, and in similar expressions below, the notation $ \Theta_{\psi}^0\bigl( \kappa \| \nabla u \|_{\LL^{A}( \Omega )}^{+} \bigr)$ stands for the one-sided limit of $\Theta_{\psi}^0(r)$ as  $r \to ( \kappa \| \nabla u \|_{\LL^{A}( \Omega )})^{+}$.
%
\end{theorem}

Let us note that the constant $c$ in inequality \eqref{measdistort} is strictly positive under assumption \eqref{scalingassump1}.
This is a consequence of the fact, by property \eqref{feb50} and Lemma \ref{thetah}, $ \Theta_{\psi}^0(r)>0$ if $r>0$

\smallskip\par 
A counterpart of inequality \eqref{measdistort} for Hausdorff contents exhibits a nonlinear dependence on the latter via   
the function $\Xi_{\psi, B} \colon (0, \infty) \to [0, \infty)$ defined as 
\begin{equation}\label{Xi}
\Xi_{\psi, B} (r) = r\, \Theta_\psi \big(\Theta_{B^{-1}} \big(1/r)^+\big) \quad \text{for $r>0$.}
\end{equation}
Observe that, as will be shown Lemma \ref{Ximon}, the function $\Xi_{\psi, B}$ is non-decreasing.
However,  it need  neither  tend to $0$ as $r \to 0^+$, nor to $\infty$ as $r\to \infty$  -- see Example \ref{ex2}.

\begin{theorem}\label{thm:contdist}
Let $n$, $m \in \N$. Let $A$ be a finite-valued  Young function satisfying conditions \eqref{feb1}, \eqref{0inf}, \eqref{feb7},
and let $B$ be the function defined by \eqref{B}.   Assume that $\varphi$ is a gauge
function fulfilling condition  \eqref{feb8}, and 
let $\psi$ be the gauge function defined by \eqref{psi}.  Let
$u \in \mathrm{V}^{1,A}(\Rn, \Rm )$.
Then 
\begin{equation}\label{Feb54}
\H^{\psi}_{\infty} \bigl( u(E) \bigr) \leq 2\kappa\Theta_{\psi}\bigl( \kappa \| \nabla u\|_{\LL^{A}(\Rn )} ^{+}\bigr)
\Xi_{\psi, B}\bigg( \frac{c_{n}\H^{\vp}_{\infty}(E)}{\kappa}\bigg)
\end{equation}
 for every set $E \subset \Rn$, where $\kappa$ is the constant appearing in   inequality \eqref{poincare}, and $c_n$ is the constant appearing in \eqref{feb25}.
\end{theorem}

\begin{example}
In the case $u \in \VV^{1,p}( \Rn , \Rm )$ for $p>n$ and $\vp (t) = t^\alpha$ for $\alpha \in (0,n]$ we have that
\begin{equation}\label{quantkauf}
  \H^{\frac{\alpha p}{\alpha +p-n}}_{\infty}\bigl( u(E) \bigr) \leq 2c_n p^{\frac{\alpha}{\alpha + p-n}}\left( \frac{p-1}{n^{\prime}-p^{\prime}} \right)^{\frac{\alpha (p-1)}{\alpha +p-n}} 
  \| \nabla u\|_{p}^{\frac{\alpha p}{\alpha +p-n}}\H_{\infty}^{\alpha}(E)^{\frac{p-n}{\alpha +p-n}}.
\end{equation}
\end{example}

Although Theorem \ref{thm:contdist} as stated  applies to the case $\vp (t) = t^n$ too, a sharper and localized result holds in this case.
This is the content of the next theorem that is stated in terms of the auxiliary function $\Phi_{B} \colon [0, \infty) \to [0, \infty)$ given by
\begin{equation}\label{Phi}
\Phi_{B} (r) = \begin{cases} r B^{-1}(1/r)^n & \quad \text{if $r>0$}
\\ 0 & \quad \text{if $r=0$.}
\end{cases}
\end{equation}
Here, $B$ denotes the function defined  in terms of $A$ and $n$ in \eqref{B}.  Under our assumptions on $A$, the function $\Phi_{B} $ is increasing, continuous 
in $[0, \infty)$ and $\lim _{r\to 0^+} \Phi_{B} (r) = 0$ -- see the proof of Lemma \ref{lemma1}. 

\begin{theorem}\label{nlusin}
Let $n, m \in \N$, with  $m \geq n$, and let $\Omega$ be an open subset of $\Rn$.
Let $A$ be a  Young function satisfying   conditions \eqref{feb1}, \eqref{0inf}, \eqref{feb7},  and let $B$ be the function defined by \eqref{B}. 
Let
$u \in \mathrm{V}^{1,A}(\Omega, \Rm )$.
Then
\begin{equation}\label{Feb60-t37}
\H^{n}_{\infty}\bigl( u(E) \bigr) \leq 2\kappa^{n+1} \| \nabla u\|_{\LL^{A}( \Omega )}^{n} \Phi_{B}\big( \Le^{n}(E)n^{\frac{n}{2}}/\kappa\big)
\end{equation}
for every set $E \subset \Omega$.
Here, $\kappa$ denotes the constant from   inequality \eqref{poincare}.
\end{theorem}

Since $\Phi_{B}\big( \Le^{n}(E)n^{\frac{n}{2}}/\kappa\big)\to 0^+$ as $\Le^{n}(E)\to 0^+$,
Theorem \ref{nlusin} in particular recaptures a stability version of the Lusin $N$ property for maps taking values in a target space
whose dimension is at least that of the domain space. Namely, it tells us that for each $\varepsilon > 0$ there exists a $\delta > 0$ such that, for sets $E \subset \Omega$,
\begin{equation}\label{generalNlusin}
\text{if \,\, $\Le^n (E) < \delta$,\,\, then \,\, $ \H^{n}_{\infty} (u(E)) < \varepsilon$.}
\end{equation}
Let us notice that inequality \eqref{Feb60-t37} is also related to the coarea formula in Sobolev spaces \cite{MSZ1}.

\medskip
 
The remaining part of this section is devoted to applications of Theorem \ref{main} to specific choices of the Young  function $A$ and of the
gauge function $\varphi$.

 \begin{example}\label{ex1}
{\rm Assume that
\begin{equation}\label{jan14}
A(t) \,\, \text{is equivalent to} \,\  t^p (\log t)^q \qquad \text{ near infinity,}
\end{equation}
where $p>n$ and $q\in \R$, and let 
\begin{equation}\label{jan15}
\varphi (r) \sim r^\alpha \big(\log \tfrac 1r)^{\beta}  \qquad \text{ near zero,}
\end{equation}
where $\alpha$ and $\beta$ fulfill one of the following alternatives 
\begin{equation}\label{jan16}
\begin{cases}
0<\alpha <n\,\, \text{and} \,\,  \beta \in \R
\\ 
\alpha =0 \,\, \text{and} \,\,   \beta <0
\\ 
 \alpha =n    \,\, \text{and} \,\,     \beta \geq  0.
\end{cases}
\end{equation}
Then property \eqref{feb6} holds with
\begin{equation}\label{jan17}
\psi(r) \sim \begin{cases} r^{\frac{\alpha p}{p+\alpha -n}} \big(\log \tfrac 1r\big)^{\frac {\alpha (q-\beta)}{p+\alpha -n}+\beta } & \quad \text{if $0\leq \alpha <n$}
\\ r^n \big(\log \tfrac 1r\big)^{\frac {\beta (p-n)}p} \big(\log (\log \tfrac 1r)\big)^{\frac {qn}p} & \quad \text{if $\alpha =n$ and $\beta >0$}
\\ r^n & \quad \text{if $\alpha =n$ and $\beta =0$}.
\end{cases} \quad \text{near zero.}
\end{equation}
In particular, if $q=0$ and $\beta=0$, then 
\begin{equation}\label{kr0}A(t)=t^p, \quad \varphi (r)=r^\alpha, \quad \psi (r) = r^{\frac{\alpha p}{p+\alpha -n}},\end{equation}
and we recover Kaufman's result  \cite[Theorem 1]{Kau}.
\\ Furthermore, if $q=0$, $\alpha =0$, then
\begin{equation}\label{sh1}
A(t)=t^p, \quad \varphi (r)\sim\big(\log \tfrac 1r\big)^{\beta}, \quad \psi (r)  \sim\big(\log \tfrac 1r\big)^{\beta}.
\end{equation}
This is an example where 
 property \eqref{Feb1} fails,    assumption  \eqref{scalingassump}  not being satisfied. More instances in the same vein are exhibited in Example \ref{ex4} below.
%
\\ In order to verify equation \eqref{jan17}, observe that,
since $i(A)=p$, condition \eqref{indexcond} is fulfilled.
Therefore,  
$$
B(t)  \,\, \text{is equivalent to} \,\  t^p (\log t)^q  \qquad \text{ near infinity.}
$$
Hence,
$$
B^{-1}(t) \sim t^{\frac 1p} (\log t)^{-\frac qp}  \qquad \text{ near infinity,}
$$
and 
$$
J(s) \sim \begin{cases} s^{\frac{p+\alpha -n}{p}} \big(\log \tfrac 1s\big)^{\frac {\beta -q}p} & \quad \text{if $0\leq \alpha <n$}
\\ s \big(\log \tfrac 1s\big)^{\frac \beta p} \big(\log (\log \tfrac 1s)\big)^{-\frac qp} & \quad \text{if $\alpha =n$ and $\beta>0$}
\\ s  & \quad \text{if $\alpha =n$ and $\beta =0$}
\end{cases} \quad \text{near zero.}
$$
Consequently,
$$
J^{-1}(r) \sim \begin{cases} r^{\frac{p}{p+\alpha -n}} \big(\log \tfrac 1r\big)^{\frac {q-\beta}{p+\alpha -n}} & \quad \text{if $0\leq \alpha <n$}
\\ r \big(\log \tfrac 1r\big)^{-\frac \beta p} \big(\log (\log \tfrac 1r)\big)^{\frac qp} & \quad \text{if $\alpha =n$ and $\beta>0$}
 \\ r   & \quad \text{if $\alpha =n$ and $\beta =0$}
\end{cases} \quad \text{near zero,}
$$
whence equation \eqref{jan17} follows.
}

\end{example}

 \begin{example}\label{ex2}
{\rm Assume that
$$
A(t) \,\, \text{is equivalent to} \,\ t^n (\log t)^q \qquad \text{ near infinity,}
$$
where $q>n-1$, and let $\varphi$ be as in \eqref{jan15}.
\\ Then property \eqref{ours} holds with 
\begin{equation}\label{pd1}
\psi(r) \sim 
\begin{cases} 
r^{\frac {n\beta}{\beta +(n-1)-q}} & \quad \text{if $ \alpha=0$}
\\ r^{n} \big(\log \tfrac 1r\big)^{q-(n-1)}  & \quad \text{if $0<\alpha < n$}
\\ r^n \big(\log (\log \tfrac 1r)\big)^{q-(n-1)}  & \quad \text{if $\alpha = n$ and $\beta >0$}
\\ r^n  & \quad \text{if $\alpha = n$ and $\beta =0$}
\end{cases} \quad \text{near zero.}
\end{equation}
We emphasize that the result concerning the case when $0<\alpha < n$ enhances 
a~result from~\cite[Theorem~1~(i)]{Raj}, where the weaker conclusion with the gauge function
$r^{n} \big(\log \tfrac 1r\big)^{q-(n-1)}$ replaced by any gauge function  $r^{n} \big(\log \tfrac 1r\big)^{\gamma}$, with $\gamma  <q-n+1$, was established.
\\ In order to prove equation \eqref{pd1}, notice that
$$
\widetilde A(t) \,\, \text{is equivalent to} \,\  t^{n'} (\log t)^{-\frac q{n-1}} \qquad \text{ near infinity,}
$$
whence
$$
\widetilde B(t)  \,\, \text{is equivalent to} \,\  t^{n'} (\log t)^{1-\frac q{n-1} } \qquad \text{ near infinity.}
$$
Therefore
$$
B(t)  \,\, \text{is equivalent to} \,\  t^{n} (\log t)^{q-(n-1)}  \qquad \text{ near infinity,}
$$
and 
\begin{equation}\label{Feb20bis}
B^{-1}(t) \sim  t^{\frac 1n} (\log t)^{\frac {(n-1)-q}n}  \qquad \text{ near infinity.}
\end{equation}
Hence,
$$
J(s) \sim \begin{cases} s^{\frac{\alpha}{n}} \big(\log \tfrac 1s\big)^{\frac {\beta+(n-1) -q}n} & \quad \text{if $0\leq \alpha <n$}
\\ s \big(\log \tfrac 1s\big)^{\frac \beta n} \big(\log (\log \tfrac 1s)\big)^{\frac {(n-1)-q}n} & \quad \text{if $\alpha =n$ and $\beta >0$}
\\ s & \quad \text{if $\alpha = n$ and $\beta =0$}
\end{cases} \quad \text{near infinity,}
$$
and
\begin{equation}\label{jan20}
J^{-1}(r) \sim \begin{cases} e^{- r^{\frac n{\beta + (n-1) -q}}}& \quad \text{if $ \alpha =0$}
\\
r^{\frac{n}{\alpha}} \big(\log \tfrac 1r\big)^{\frac {q-(n-1)-\beta}{\alpha}} & \quad \text{if $0< \alpha <n$}
\\ r \big(\log \tfrac 1r\big)^{-\frac \beta n} \big(\log (\log \tfrac 1r)\big)^{\frac {q-(n-1)}n} & \quad \text{if $\alpha =n$ and $\beta >0$ }
\\ r & \quad \text{if $\alpha = n$ and $\beta =0$}
\end{cases} \quad \text{near zero.}
\end{equation}
Equation \eqref{pd1} follows from \eqref{jan15} and \eqref{jan20}.
\\
Let us notice that, if
either $0<\alpha <n$ and $\beta \in \R$, or $\alpha =n$ and $\beta >0$, then the function given by \eqref{Xi} does not tend to $0$ as $r \to 0^+$.
Indeed, 
 the function $\varphi$ fulfills condition \eqref{Feb4}. Moreover, by equation \eqref{Feb20bis},
$$\Theta_{B^{-1}}(1/r) \geq \Theta_{B^{-1}}^\infty (1/r) = r^{-\frac 1n} \quad \text{for $r>0$,}$$
and, by equation \eqref{pd1},
$$\Theta_{\psi} (s) \geq \Theta_{\psi}^0 (s) = s^n  \quad \text{for $s>0$.}$$
Hence, by the monotonicity of the functions $\Theta_{B^{-1}}$ and $\Theta_{\psi}$,
$$ \Xi_{\psi, B}(r)= r\, \Theta_{\psi}\left( \Theta_{B^{-1}} (1/r)^+ \right) \geq 1 \quad  \text{for $r>0$.}$$
}
\end{example}

 \begin{example}\label{ex3}
{\rm Assume that
$$
A(t) \,\, \text{is equivalent to} \,\, e^{t^\gamma}\qquad \text{ near infinity,}
$$
for some  $\gamma >0$, and let $\varphi$ be as in \eqref{jan15}.
Then equation \eqref{feb6} holds with
\begin{equation}\label{jan21}
\psi(r) \sim 
\begin{cases} 
 r^{\alpha} \big(\log \tfrac 1r\big)^{\beta - \frac \alpha \gamma}   & \quad \text{if $ 0\leq \alpha <n$}
\\ r^{n} \big(\log \tfrac 1r\big)^{\beta} \big(\log (\log \tfrac 1r)\big)^{-\frac {n}\gamma} & \quad \text{if $\alpha= n$ and $\beta >0$}
\\  r^n & \quad \text{if $\alpha= n$ and $\beta =0$}
\end{cases} \quad \text{near zero.}
\end{equation}
Indeed, since $i(A)=\infty$, condition \eqref{indexcond} holds, whence
$$
B(t)  \,\, \text{is equivalent to} \,\,   e^{t^\gamma}  \qquad \text{ near infinity.}
$$
Hence, 
$$
B^{-1}(t) \sim  (\log t)^{\frac {1}\gamma}  \qquad \text{ near infinity.}
$$
Thus,
$$
J(s) \sim \begin{cases} s(\log \tfrac 1s)^{\frac {1}\gamma} & \quad \text{if $0\leq \alpha <n$}
\\ s   \big(\log (\log \tfrac 1s)\big)^{\frac {1}\gamma} & \quad \text{if $\alpha =n$ and $\beta >0$}
\\ s & \quad \text{if $\alpha =n$ and $\beta =0$}
\end{cases} \quad \text{near infinity,}
$$
and
\begin{equation}\label{jan22}
J^{-1}(r) \sim \begin{cases}  
r(\log \tfrac 1r)^{-\frac {1}\gamma} & \quad \text{if $0\leq \alpha <n$}
\\ r   \big(\log (\log \tfrac 1r)\big)^{-\frac {1}\gamma} & \quad \text{if $\alpha =n$ and $\beta >0$}
\\ r    & \quad \text{if $\alpha =n$ and $\beta =0$}
\end{cases} \quad \text{near zero.}
\end{equation}
Coupling equation \eqref{jan15} with \eqref{jan22} yields \eqref{jan21}.}

\end{example}

\begin{example}\label{ex4}
{\rm Here, we exhibit classes of Young functions $A$ and gauge functions $\varphi$ for which the gauge function $\psi$
appearing in property \eqref{feb6} is equivalent to $\varphi$. Loosely speaking, this is the case when the function $A$ grows
essentially faster than the $n$-th power  near infinity, and the function $\varphi$ decays more slowly than any power near zero.
These instances provide us with a generalization of example \eqref{sh1}.
Clearly, they also show that property \eqref{Feb1} may actually fail if assumption  \eqref{scalingassump} is not fulfilled.
\\ Precisely, assume that $A$ is any Young function such that 
$$
i_\infty(A)>n,
$$
and let $\varphi$ be any gauge function such that
\begin{equation}\label{jan40'}
\liminf_{r \to 0^+} \frac{ \varphi (r^\gamma)}{ \varphi (r)} > 0  \quad \text{for every $\gamma >0$.}
\end{equation}
Then property   \eqref{feb6} holds with
\begin{equation}\label{jan41'}
\psi(r) \sim \varphi  (r)\quad \text{near zero.}
\end{equation} 
To verify this assertion, recall that the current assumption on $A$ ensures that $B$ is equivalent to $A$ near infinity.
Moreover, it also implies that there exists $\epsilon >0$ such that $A(t) \geq t^{n+\epsilon} $ near infinity. Altogether, one has that
\begin{equation}\label{jan43}
B(t) \geq t^{n+\epsilon} \quad \text{near infinity.}
\end{equation}
Now, in order to establish equation \eqref{jan41'} it suffices to show that
\begin{equation}\label{jan46}
\psi(r) \geq  c \,\varphi  (r) \quad \text{near zero}
\end{equation} 
for some positive constant $c$,
the reverse inequality being always fulfilled (up to a multiplicative constant). Owing to assumption \eqref{jan40'},
inequality \eqref{jan46} will in turn follow if we prove that there exists $\gamma >0$ such that
\begin{equation}\label{jan47}
J^{-1}(r) \geq r^\gamma \quad \text{near zero,}
\end{equation} 
or, equivalently, 
\begin{equation}\label{jan48}
J(s) \leq s^{\frac 1\gamma} \quad \text{near zero.}
\end{equation} 
The last inequality reads
\begin{equation}\label{jan49}
 s\, B^{-1}\Big(\frac {\varphi (s)}{s^n}\Big)\leq s^{\frac 1\gamma}  \quad   \text{near zero,}
\end{equation}
namely
\begin{equation}\label{jan49--}
\varphi (s) \leq s^n B\big(s^{\frac 1\gamma-1}\big)  \quad   \text{near zero.}
\end{equation}
Choosing $\gamma > {\color{blue} 1 + \frac{n}\epsilon}$ and making use of equation \eqref{jan43} tell us that
\begin{equation}\label{jan50}
 s^n B\big(s^{\frac 1\gamma-1}\big) \geq s^{\frac{n+\epsilon}\gamma -\epsilon} \quad   \text{near zero.}
\end{equation}
By our choice of $\gamma$, one has that $\frac{n+\epsilon}\gamma -\epsilon <0$. Inequality \eqref{jan49--} hence follows.
}

\end{example}

\section{Technical lemmas}\label{tech}

The results collected in this section concern properties of the functions $B$, $J$ and $\psi$, defined by equations  \eqref{B}, \eqref{J} and \eqref{psi},   of the scaling functions given by \eqref{defthetah}-- \eqref{defthetahs}, and of the function $\Xi_{\psi, B}$ defined as in \eqref{Xi}.

The first result summarizes certain basic properties of  the function  $A$ which follow from our assumptions, and tells us that they are  inherited by $\widetilde A$ and  by $B$.

\begin{lemma}\label{Bcont}
Assume that $A$ is a   Young function fulfilling conditions \eqref{feb1}, \eqref{0inf} and \eqref{feb7}. 
\\ (i) We have that
\begin{equation}\label{aug200}
\lim_{t\to 0^+} \frac{A(t)}t =0 \qquad \text{and} \qquad \lim_{t\to \infty} \frac{A(t)}t =\infty.
\end{equation}
Moreover, the function
\begin{equation}\label{aug201}
\text{ $\frac{A(t)}t $ \,\, is increasing in $(0, \infty)$.} 
\end{equation}
\\ (ii) The Young conjugate $\widetilde A$ is such that $0<\widetilde A (t) <\infty$  for $t>0$, and 
\begin{equation}\label{aug202}
 \lim_{t\to 0^+} \frac{\widetilde A(t)}t =0 \quad \text{and} \qquad \lim_{t\to \infty} \frac{\widetilde A(t)}t =\infty.
\end{equation}
\\ (iii) 
The Young function $B$  given by \eqref{B} is such that  $0<B (t) <\infty$  for $t>0$, and 
\begin{equation}\label{aug203}
 \lim_{t\to 0^+} \frac{B(t)}t =0 \quad \text{and} \qquad \lim_{t\to \infty} \frac{B(t)}t =\infty.
\end{equation}
In particular, the function $B: [0, \infty) \to [0, \infty)$ is bijective. 
%
%
\end{lemma}
\textit{Proof} (i) The first limit in equation \eqref{aug200} follows from assumption \eqref{feb7} and property \eqref{monotone}. The second limit is a consequence of the same property and of assumption  \eqref{feb1}. The strengthening of property \eqref{monotone} stated in equation \eqref{aug201} holds owing to the first limit in \eqref{aug200}.
\\ (ii) The properties of the function $\widetilde A$ claimed here are classically known to follow, via its very definition, from parallel properties enjoyed by $A$.
\\ (iii)
The statement is trivial when $n=1$, thanks to Part (i). Assume that $n \geq 2$. As observed in the proof of Part (ii), it suffices to prove the statement with $B$ replaced by $\widetilde B$. Since
$$\widetilde B(t) =  t^{n'}\displaystyle \int_t^\infty \frac{\widetilde A(s)}{s^{1+n'}}\, ds \quad \text{for $t>0$,}$$
one clearly has that $0<\widetilde B(t)<\infty$ for $t>0$. Moreover, thanks to equation \eqref{aug202},
$$\lim_{t \to 0^+} \frac{\widetilde B(t)}t = \lim_{t \to 0^+}  t^{n'-1}\displaystyle \int_t^\infty \frac{\widetilde A(s)}{s^{1+n'}}\, ds =0,$$
and
$$\lim_{t \to \infty} \frac{\widetilde B(t)}t = \lim_{t \to \infty}  t^{n'-1}\displaystyle \int_t^\infty \frac{\widetilde A(s)}{s^{1+n'}}\, ds =\infty.$$
\qed

The following result ensures that the inverse $J^{-1}$ of $J$ is well defined.

\begin{lemma}\label{lemma1}
Let $n \in \N$,
let $A$ be a Young function fulfilling conditions \eqref{feb1}, \eqref{0inf}, \eqref{feb7}, and let $\varphi$ be a gauge function.
Then the function $J\colon (0, \infty) \to (0, \infty)$, given by \eqref{J},
is increasing and such that
\begin{equation}\label{feb9}
\lim_{t\to 0^+}J(t)=0 \qquad \text{and} \qquad \lim_{t\to \infty}J(t)=\infty. 
\end{equation}
In particular, $J$ is bijective.
\end{lemma}
\textit{Proof} 
To begin with, by Lemma \ref{Bcont} the Young function  $B$ is bijective, and hence its inverse $B^{-1} : [0, \infty) \to [0, \infty)$ is classically well defined.
Suppose first that $n \geq 2$. One has that the function
\begin{equation}\label{feb10}
\frac{B(t)}{t^n} \quad \text{is non-decreasing}
\end{equation}
Indeed, property \eqref{feb10} is equivalent to the fact that the function
\begin{equation}\label{feb11}
\frac{\widetilde B(t)}{t^{n'}}  \quad \text{is non-increasing}
\end{equation}
and the latter holds, since
\begin{equation}\label{feb12}
\frac{\widetilde B(t)}{t^{n'}}= \int_t^\infty \frac{\widetilde A(s)}{s^{1+n'}}\, ds \qquad \text{for $t>0$.}
\end{equation}
As a consequence of property \eqref{feb10}, 
the function
\begin{equation}\label{feb13}
tB^{-1}(t^{-n})  \quad \text{is non-decreasing.}
\end{equation}
Hence, the function
\begin{equation}\label{feb14}
  J(t)=  tB^{-1}\Big(\frac{\varphi (t)}{t^n}\Big) = \frac{t}{\varphi (t)^{\frac 1n}} B^{-1}\Big(\frac{\varphi (t)}{t^n}\Big)
  \varphi (t)^{\frac 1n}  \quad \text{is increasing,}
\end{equation}
inasmuch as the function $\varphi$ is increasing and satisfies condition \eqref{feb8}.
\\ Next, owing to equations \eqref{feb3} and \eqref{feb14}, and to the fact that the function $\varphi$ is increasing, we have that
\begin{equation}\label{feb60}
\lim_{t\to 0^+} J(t) \leq \lim_{t\to 0^+} t B^{-1}\Big(\frac{1}{t^n}\Big).
\end{equation}
Moreover, a property of Young functions tells us that
\begin{equation}\label{feb144}
t \leq B^{-1} (t) \widetilde B^{-1}(t) \leq 2t \qquad \text{for $t \geq 0$.}
\end{equation}
Thus,
$$\lim _{t\to 0^+} J(t)=0 \quad \text{if} \quad \lim_{t\to \infty} \frac{B^{-1}(t)}{t^{\frac 1n}} = 0.$$
Now, notice that
\begin{align}\label{feb15}
\lim_{t\to \infty} \frac{B^{-1}(t)}{t^{\frac 1n}} = 0 
 \quad \text{if and only if} \quad \lim_{t \to \infty} \frac{t^{\frac 1{n'}}}{\widetilde B^{-1}(t)}=0
\quad 
 \text{ if} \quad \lim_{t \to \infty}\frac{\widetilde B(t)}{t^{n'}}=0,
\end{align}
and the last limit holds, owing to equation \eqref{feb12}. The first limit in \eqref{feb9} is thus established. 
\\ As for the second one,  from equations \eqref{feb3} and \eqref{feb14} again, and   the   monotonicity of the function $\varphi$,  one has that
\begin{equation}\label{feb61}
\lim_{t\to \infty} J(t) \geq {\color{blue}c}\lim_{t\to \infty} tB^{-1}\Big(\frac{1}{t^n}\Big),
\end{equation}
where $c=\min\bigl\{1,\lim\limits_{t\to\infty}\varphi(t)\bigr\}$.  Therefore, by an analogous chain as above,
\begin{equation}\label{jan11}
\lim _{t\to \infty} J(t) = \infty \quad \text{if} \quad  \lim_{t\to 0^+} \frac{B^{-1}(t)}{t^{\frac 1n}} = \infty 
\quad \text{if and only if} \quad    \lim_{t\to 0^+} \frac{t^{\frac 1{n'}}}{\widetilde B^{-1}(t)} = \infty
\quad \text{if} \quad
     \lim_{t \to 0^+}\frac{\widetilde B(t)}{t^{n'}}=\infty.
\end{equation}
By \cite[Lemma 4.4]{Cianchi_Ibero}, condition \eqref{feb7} is equivalent to
\begin{equation}\label{feb7bis}
\int_0 \frac{\widetilde A(t)}{t^{1+n'}}\, dt = \infty.
\end{equation}
Hence, the last limit in chain \eqref{jan11} holds. This proves the second limit in \eqref{feb9}.
\\
These pieces of information conclude the proof,
since the function $J$ is continuous,  inasmuch as the function $B^{-1}$ is continuous in $(0, \infty)$.
\\ Assume next that $n=1$, and hence $B=A$. Since $A$ is a Young function,  equation \eqref{feb13} still holds. Hence, we conclude as above that $J$ is increasing. The limits 
$$\lim _{t\to 0^+} J(t)=0  \quad \text{and} \quad \lim_{t\to \infty} J(t) =\infty$$
follow from equations \eqref{feb60} and \eqref{feb61}, with $B=A$, via property \eqref{aug200}.
\qed

The next lemma substantiates the fact that $\psi$, defined as in \eqref{psi}, is actually a gauge function.

\begin{lemma}\label{lemma3}
Let $n$, $A$ and $\varphi$ be as in Lemma \ref{lemma1}. Then  the function $\psi$ defined by \eqref{psi} is continuous
and increasing, and the function
\begin{equation}\label{feb50}
\text{$\frac{\psi (r)}{r^n}$ is non-increasing.}
\end{equation}
In particular, 
\begin{equation}\label{delta2psi}
\psi (k r) \leq \max\{1, k^n\} \psi (r) \qquad \text{for $k>0$ and $r \geq 0$.}
\end{equation}
\end{lemma}
\textit{Proof} 
The fact that  $\psi$ is continuous and increasing is a consequence of the fact that both $\varphi$ and
$J$ are continuous and increasing. As for property \eqref{feb50}, note that it is equivalent to the fact that the function
\begin{equation}\label{feb51}
\text{ $\frac{\varphi (r)}{J(r)^n}$ is non-increasing.}
\end{equation}
On the other hand,
\begin{equation}\label{feb52}
\frac{\varphi (r)}{J(r)^n} = \frac{\frac{\varphi (r)}{r^n}}{B^{-1}\big(\frac{\varphi (r)}{r^n}\big)} \qquad \text{for $r>0$.}
\end{equation}
By property \eqref{feb8} the function $\frac{\varphi (r)}{r^n}$ is  non-increasing, whereas the function $\frac t{B^{-1}(t)}$
is non-decreasing, inasmuch as $B$ is a Young function (in fact, it is increasing, as a consequence of Lemma \ref{Bcont}, Part (iii)). Hence, property \eqref{feb51} follows.
\\ Inequality \eqref{delta2psi} is a consequence of inequality \eqref{delta2}, which, thanks to   property \eqref{feb50},  can be applied with $\varphi$ replaced by $\psi$.
\qed

The   function $\psi$ owes its definition to the inequality established in the following lemma.

\begin{lemma}\label{lemma2}
Let $n$, $A$ and $\varphi$ be as in Lemma \ref{lemma1}, and let $B$ and $\psi$ be the functions defined by \eqref{B} and \eqref{psi}, respectively.  Then 
\begin{equation}\label{feb16}
\psi (st) \leq  \varphi (t) +  t^n B(s)  \qquad \text{for $s, t>0$.}
\end{equation}
\end{lemma}
\textit{Proof} Inequality \eqref{feb16} trivially holds if either $s=0$ or $t=0$. On the other hand,  if both $s>0$ and  $t>0$, it is equivalent to 
\begin{equation}\label{feb17}
\psi (s) \leq  \varphi (t) +  t^n B(s/t).
\end{equation}
The function $\varphi$ is increasing, by our assumption. On the other hand, by property \eqref{feb10}, 
for each $s>0$, the function $t^n B(s/t)$ is non-increasing in $t$. Now
  observe that
\begin{equation}\label{feb18}
\varphi (J^{-1} (s)) = J^{-1} (s)^n B(s/J^{-1} (s)) \quad \text{for $s>0$.}
\end{equation}
The monotonicity properties mentioned above imply that
$$ \varphi (t) +  t^n B(s/t) \geq \varphi (J^{-1} (s)) \qquad \text{if $t \geq J^{-1} (s)$,}$$
and 
$$ \varphi (t) +  t^n B(s/t) \geq  J^{-1} (s)^n B(s/J^{-1} (s))  \qquad  \text{if $0 < t \leq J^{-1} (s)$.}$$
Hence, by equation \eqref{feb18} and the definition of the function $\psi$,
$$ \varphi (t) +  t^n B(s/t) \geq \varphi (J^{-1} (s)) = \psi (s)  \qquad \text{for $t>0$.}$$
Inequality \eqref{feb17} is thus established.
\qed

Our analysis of the scaling functions begins with the next result, which, in particular, deals with a submultiplicativity property.

\begin{lemma}\label{thetah}
Let $h \colon (0,\infty ) \to (0,\infty )$ be a non-decreasing continuous function
 such that  the function 
\begin{equation}\label{Feb36}
\text{$\frac{h(t)}{t^\gamma}$ is non-increasing  for some $\gamma >0$.}
\end{equation}
Denote by $\Theta$ either of the functions $ \Theta_{h}^0$,  $\Theta_h^\infty$ and   $\Theta_{h}$ associated with $h$ as in \eqref{defthetah}, \eqref{defthetaH} and \eqref{defthetahs}, respectively, and by $\Theta_*$ the function defined analogously, with $\lq\lq \limsup"$ or $\lq\lq \sup"$ replaced with $\lq\lq \liminf"$ or $\lq\lq \inf"$, respectively.  Then
\begin{equation}\label{Feb30}
\min \{ 1, r^\gamma \} \leq \Theta (r) \leq \max \{ 1,r^\gamma \} \quad \text{for $r>0$,}
\end{equation}
and
\begin{equation}\label{Feb32}
\Theta (rs) \leq \Theta (r) \Theta (s) \quad \text{for $r,s>0$.}
\end{equation}
 Furthermore,
\begin{equation}\label{Feb33}
\Theta_* (r) = \frac{1}{\Theta \bigl( \tfrac{1}{r}\bigr)}  \quad \text{for $r>0$.}
\end{equation}
\end{lemma}

\medskip

\noindent
\textit{Proof} We consider the case when $\Theta = \Theta_h^0$, the proof in the remaining cases being  analogous. The monotonicity of the function $h$ and property \eqref{Feb36} 
 ensure that $r^\gamma \leq \Theta_{h}^0(r)  \leq 1$ for $r \in (0,1]$.
Indeed, if $r \in (0,1]$, then the upper bound follows since $h$ is non-decreasing, whereas  the lower
bound is a consequence of the fact that
$$
\frac{h(rt)}{h(t)} = \frac{h(rt)}{(rt)^\gamma} \, \frac{t^\gamma}{h(t)} r^\gamma \geq r^\gamma
$$
for $t>0$, whence $\Theta_{h}^0(r) \geq r^\gamma$. The proof of the inequalities $1 \leq \Theta_{h}^0(r)  \leq r^\gamma$ for $r > 1$ is analogous. Equation \eqref{Feb30} is therefore established.
\\ As for inequality \eqref{Feb32}, 
 fix $r>0$ and $s>0$, and choose a sequence $\{t_k\}$ such that $t_k \to 0^+$ and $h(rst_{k})/h(t_k ) \to \Theta_{h}^0(rs)$.
Now
$$
\Theta_{h}^0(rs) = \lim_{k \to \infty} \frac{h(rst_{k} )}{h(st_{k})} \, \frac{h(st_{k})}{h(t_{k})}
\leq  \Theta_{h}^0(r)\Theta_{h}^0(s),
$$
namely \eqref{Feb32}. 
\\ Finally, equation \eqref{Feb33} follows from  the  chain:
$$
\frac{1}{\Theta_{h}^0\left( \tfrac{1}{r}\right)} = \left( \limsup_{t \to 0^{+}}\frac{h(t/r)}{h(t)} \right)^{-1}
= \liminf_{t \to 0^{+}} \frac{h(t)}{h(t/r)} = \liminf_{t \to 0^{+}} \frac{h(rt)}{h(t)} = \big(\Theta^0_h\big)_*(r) \quad \text{for $r>0$.}
$$  \qed

\begin{corollary}\label{dicotomy}
Let $h \colon (0,\infty ) \to (0,\infty )$ be a non-decreasing continuous function and let $\Theta$ be any of  the functions associated with $h$ as in \eqref{defthetah}--\eqref{defthetahs}. Then
\begin{equation}\label{dicot}
\text{either \quad  $\lim_{r \to 0^+} \Theta(r) =0$, \quad  or \quad $\Theta(r)=1$ \quad  for every $r\in (0, 1]$.}
\end{equation}
\end{corollary}

\medskip

\noindent
\textit{Proof}
One  
has that  $\Theta (r) \in (0,1]$ for  $r \in (0,1]$. If there exists $r_{0} \in (0,1)$ such that $\Theta(r_{0}) < 1$, then, by property \eqref{Feb32},    $\Theta (r_{0}^{k}) \leq \Theta (r_{0})^{k}$ for every $k \in \N$. Hence,   $\lim_{r \to 0^+} \Theta(r)=0$.
\qed

The next two lemmas are concerned with properties of a one-parameter family of functions $J_r$ obtained by inserting a scaling factor $r$ in the definition of function $J$. Specifically, 
 the function $J_{r} \colon [0,\infty) \to [0,\infty)$ is defined, for $r>0$, as
\begin{equation}\label{Jr}
J_{r}(s) = \begin{cases} sB^{-1}\left(\displaystyle \frac{r\vp (s)}{s^n} \right) & \quad \text{for $s > 0$}
\\ 0  & \quad \text{for $s = 0$.}
\end{cases}
\end{equation}

\begin{lemma}\label{aux1}
Let $n \in \N$,
let $A$ be a Young function fulfilling conditions \eqref{feb1}, \eqref{0inf} and \eqref{feb7}, let $B$ be the function given by \eqref{B}, and let $\varphi$ be a gauge function. 
Then, for each fixed $r$, the function $J_r$ given by \eqref{Jr} is increasing and  bijective, and   for each fixed $s>0$ the function 
$$(0, \infty) \ni r \mapsto J_{r}^{-1}(s) \quad \text{is non-increasing.}$$
Moreover,
\begin{equation}\label{modifyfeb16}
\vp \bigl( J_{r}^{-1}(st) \bigr) \leq \vp (t) + \frac{t^n}{r}B(s)\quad \text{for $s,t >0$.}
\end{equation}
\end{lemma}
 
\noindent
\textit{Proof.}
The function $J_r$ is defined  as the function $J$ given by \eqref{J}, with the gauge function $\vp$   replaced by $r\vp$. 
The asserted properties of  the function $J_r$ follow from Lemmas \ref{lemma1} and \ref{lemma2}. In particular, inequality \eqref{modifyfeb16} follows on replacing the function $\vp$ with $r\vp$ in equation 
 \eqref{feb16},  and dividing through the resultant inequality by $r$.
\qed
 
\begin{lemma}\label{aux2}
Let $n \in \N$,
let $A$ be a Young function fulfilling conditions \eqref{feb1}, \eqref{0inf}, \eqref{feb7}, and  let $B$ be the function given by \eqref{B}. Assume that $\varphi$ is a gauge function satisfying property \eqref{Feb4}, and let $\psi$ be the gauge function defined by \eqref{psi}. Let $J_r$ be the function defined by \eqref{Jr}. 
 Then
\begin{equation}\label{FebA}
\liminf_{t \to 0^{+}}\frac{\vp \bigl( J_{r}^{-1}(t) \bigr)}{\psi (t)} \geq \frac{1}{\Theta_{\psi}^0\bigl( \Theta_{B^{-1}}^\infty(r)^{+} \bigr)}\quad \text{for  $r >0$,}
\end{equation}
and 
\begin{equation}\label{FebB}\inf_{t>0}\frac{\vp \bigl( J_{r}^{-1}(t) \bigr)}{\psi (t)} \geq \frac{1}{\Theta_{\psi}\bigl( \Theta_{B^{-1}}(r)^{+} \bigr)} \quad \text{for  $r>0$.}
\end{equation}
\end{lemma}

\noindent
\textit{Proof.} We limit ourselves to proving equation \eqref{FebA}, the proof of \eqref{FebB} being analogous. Fix $r>0$. By the very definitions of the functions $\psi$, $J$ and $J_r$, and the properties of the latter established in Lemma \ref{aux2},  one has that
$$
\liminf_{t \to 0^{+}}\frac{\vp \bigl( J_{r}^{-1}(t) \bigr)}{\psi (t)} =  \liminf_{s \to 0^{+}}\frac{\vp(s)}{\psi (J_r(s))} =
\left( \limsup_{s \to 0^{+}} \frac{\psi \left( s B^{-1}\bigl( r\tfrac{\vp (s)}{s^{n}} \bigr)\right)}{\psi \left( sB^{-1} \bigl( \tfrac{\vp (s)}{s^{n}} \bigr) \right)}\right)^{-1}.
$$
Therefore, inequality \eqref{FebA} will  follow if we  show that
\begin{equation}\label{Feb43}
\limsup_{s \to 0^z{+}} \frac{\psi \left( sB^{-1}\bigl( r\tfrac{\vp (s)}{s^{n}} \bigr)\right)}{\psi \left( sB^{-1}\bigl( \tfrac{\vp (s)}{s^{n}} \bigr)\right)} 
\leq \Theta_{\psi}^0 \bigl( \Theta_{B^{-1}}^{\color{blue}\infty}(r)^{+} \bigr).
\end{equation}
To verify inequality \eqref{Feb43}, observe that, owing to assumption \eqref{Feb4},  for every $\varepsilon > 0$ there exists $\delta > 0$ such that
$$
B^{-1}\bigl( r\tfrac{\vp (s)}{s^n} \bigr) < \Theta_{B^{-1}}^\infty (r)(1+ \varepsilon )B^{-1}\bigl( \tfrac{\vp (s)}{s^n} \bigr) \quad \text{for  $s \in (0, \delta )$.}
$$
Hence, inasmuch as  $\psi$ is an increasing function,  
$$
\frac{\psi \left( sB^{-1}\bigl( r\tfrac{\vp (s)}{s^n} \bigr) \right)}{\psi \left( sB^{-1}\bigl( \tfrac{\vp (s)}{s^n} \bigr)\right)} \leq \frac{\psi \left( \Theta_{B^{-1}}^\infty(r)(1+\varepsilon )
sB^{-1}\bigl( \tfrac{\vp (s)}{s^n} \bigr) \right)}{\psi \left( sB^{-1}\bigl( \tfrac{\vp (s)}{s^n} \bigr) \right)}  \quad \text{for  $s \in (0, \delta )$.}
$$
Since, by equation \eqref{feb9}, $\lim _{s \to 0^{+}} sB^{-1}\bigl( \tfrac{\vp (s)}{s^n} \bigr)=0$, we deduce from this inequality that 
$$
\limsup_{s \to 0^{+}} \frac{\psi \left( sB^{-1}\bigl( r\tfrac{\vp (s)}{s^{n}} \bigr)\right)}{\psi \left( sB^{-1}\bigl( \tfrac{\vp (s)}{s^{n}} \bigr)\right)} 
 \leq \Theta_{\psi}^0 \bigl( \Theta_{B^{-1}}^\infty(r)(1+ \varepsilon ) \bigr).
$$
Inequality \eqref{Feb43} hence follows, owing to the arbitrariness of  $\varepsilon $.
\qed

\begin{remark}{\rm 
Under the additional condition  that $i_\infty (A)>n$,   the Young function $B$ defined by \eqref{B} is equivalent to $A$ at infinity.  Namely,
\begin{equation}\label{Feb80}
B \bigl( c_{1}t \bigr) \leq A(t) \leq B \bigl( c_{2}t \bigr)
\end{equation}
 for every $t \geq t_0$, where $c_1$, $c_2$, $t_0$ are positive constants. As a consequence,
\begin{equation}\label{Feb20}
\Theta_{A^{-1}}^\infty(0^{+}) =0 \, \mbox{ if and only if } \, \Theta_{B^{-1}}^\infty(0^{+})=0.
\end{equation}
Indeed, inequalities \eqref{Feb80} imply that
\begin{equation}\label{Feb81}
\frac{B^{-1}(s)}{c_2} \leq A^{-1}(s) \leq \frac{B^{-1}(s)}{c_1}
\end{equation}
for large $s$. Hence,
 $\frac{c_1}{c_2}\Theta_{B^{-1}}^\infty(r)\leq \Theta_{A^{-1}}^\infty(r) \leq \frac{c_2}{c_1}\Theta_{B^{-1}}^\infty(r)$ for $r >0$, whence assertion \eqref{Feb20} follows. 
\\ Analogously, if $i(A) > n$, then 
\begin{equation}\label{Feb20'}
\Theta_{A^{-1}}(0^{+}) =0 \, \mbox{ if and only if } \, \Theta_{B^{-1}}(0^{+})=0.
\end{equation}}
\end{remark}

\medskip

We conclude with a proof of the monotonicity of the function $\Xi_{\psi, B}$.

\begin{lemma}\label{Ximon}
Let $n \in \N$,
let $A$ be a Young function fulfilling conditions \eqref{feb1}, \eqref{0inf},\eqref{feb7},  and let $B$ be the function given by \eqref{B}. Assume that $\varphi$ is a gauge function and let $\psi$ be the function given by \eqref{psi}. Then the function $\Xi_{\psi, B}$, defined by \eqref{Xi}, is non-decreasing.
\end{lemma}

\noindent
\textit{Proof.}  It clearly suffices to show that, for fixed $s, t>0$, the function
\begin{equation}\label{scale2}
(0,\infty ) \ni r \mapsto \frac{\psi \left( \frac{B^{-1}\left( \frac{t}{r} \right)}{B^{-1}(t)}s \right)}{\psi (s)}r \quad \text{is non-decreasing.}
\end{equation}
On the other hand, property  \eqref{scale2} holds if the function
$$
(0,\infty ) \ni r \mapsto\psi \left( \frac{B^{-1}(rt)}{B^{-1}(t)}s \right)\,  \frac 1r  \quad \text{is non-increasing.}
$$
 Since
$$
  \psi \left( \frac{B^{-1}(rt)}{B^{-1}(t)}s \right) \,  \frac 1r  = \frac{\psi \left( B^{-1}(rt)\displaystyle\frac{s}{B^{-1}(t)}\right)}{\left( B^{-1}(rt)\displaystyle\frac{s}{B^{-1}(t)}\right)^{n}}
\cdot \frac{\left( B^{-1}(rt)\displaystyle\frac{s}{B^{-1}(t)}\right)^{n}}{r}\quad \text{for $r>0$,}
$$
this property is a consequence of equations \eqref{feb50} and \eqref{feb13}.
\qed

\section{Proofs of the main results}\label{proofs}

Having  the necessary preliminary material at disposal, we are now in a position to accomplish the proofs of our main results.
 
\medskip
\par\noindent
\textit{Proof of Theorem \ref{main}.}   Let us begin with Part (i) and
let us temporarily assume, in addition,   that the closure $\overline E$ of $E$ is 
compact and $\overline E \subset \Omega$. Hence, there exists 
$\varrho_0 > 0$ such that the metric neighbourhood of $E$ of radius $\varrho_0$, defined as
 \begin{equation}\label{Brho}
B_{\varrho_0}(E)= \{x\in \Rn: {\rm dist }(x, E) < \varrho_0\},
\end{equation}
is compactly contained in $\Omega$. 
 Owing to our assumptions on $A$, the function $u$ is uniformly continuous in $B_{\varrho_0}(E)$. Let us denote by $\omega$  its modulus of continuity.
\\ Let  $\kappa $  be the constant appearing in   inequality \eqref{poincare} and let  $\vp_{\kappa}$ be
the function defined as in \eqref{phik}. By virtue of equations  \eqref{feb25} and \eqref{feb26}, there exists a constant $c=c(n, \kappa )$ such
that $\Lambda^{\vp_{\kappa}}(E) \leq c\H^{\vp}(E)$. Now,  suppose that
$\H^{\vp}(E) < \delta$, where the number $\delta>0$ will be specified in the course of the proof.  Hence $\Lambda^{\vp_{\kappa}}(E) < c\delta$
and for any   $\sigma \in (0,1)$ we may   select a family of non-overlapping dyadic cubes $\{ Q_i \}$ such that $E \subset \bigcup_{i} Q_i$, 
\begin{equation}\label{feb30}
\sup_{i} d(Q_i) < \sigma,
\end{equation}
and
\begin{equation}\label{feb27}
\sum_i \varphi (\kappa d(Q_i)) <  c\delta.
\end{equation}
We may clearly also  assume that $Q_i \subset B_{\varrho_0}(E)$. 
Hence,  
\begin{equation}\label{jan26}
 \sup_{i} d(u(Q_{i})) < \omega(\sigma).
\end{equation}
From inequalities \eqref{feb16} and \eqref{poincare} we infer that, if $\lambda > \|\nabla u\|_{L^A(\Omega)}$, then
\begin{align}\label{feb28}
\psi\bigg(\frac{d(u(Q_i))}\lambda \bigg) \leq \varphi (\kappa d(Q_i)) + \kappa ^n d(Q_i)^nB\bigg(\frac{d(u(Q_i))}{\kappa \lambda d(Q_i)}\bigg)
\leq  \varphi (\kappa d(Q_i)) + \kappa^{n+1} \int_{Q_i}A\bigg(\frac{|\nabla u|}\lambda\bigg)\, dx
\end{align}
for every $i$. Hence, by Lemma \ref{lemma3},
\begin{align}\label{feb29}
\min\{1, \lambda^{-n}\} \sum_i \psi\big(d(u(Q_i))\big) &  \leq \sum_i \psi\bigg(\frac{d(u(Q_i))}\lambda \bigg)  \\ \nonumber &
\leq \sum_i \varphi (\kappa d(Q_i)) + \kappa^{n+1}   \int_{\bigcup_i Q_i} A\bigg(\frac{|\nabla u|}\lambda\bigg)\, dx.
\end{align}
Owing to equations \eqref{feb8}, \eqref{feb30} and \eqref{feb27},  
\begin{align}\label{feb31}
  \Le^n \Big(\bigcup_i Q_i\Big)  &  = n^{-\frac{n}{2}}\sum_i d(Q_i)^n = n^{-\frac{n}{2}}\kappa^{-n}\sum_i \vp \big(\kappa d(Q_i)\big)
  \frac{(\kappa d(Q_i))^n}{\vp \big(\kappa d(Q_i)\big)}  \\ \nonumber & \leq \frac {c\delta}{n^{\frac{n}{2}}\kappa^n}
  \frac{\bigl( \kappa\sigma  \bigr)^{n}}{\vp \bigl( \kappa  \sigma \bigr)}
 \leq   \frac {c\delta}{n^{\frac{n}{2}}\varphi (\kappa)}.
\end{align}
Hence, there exists   $\delta' = \delta' (n,m,A, \varphi, u, \varepsilon)>0$ such that
\begin{equation}\label{feb32}
\kappa^{n+1}\int_{\bigcup_{i}Q_i} \! A\left( \frac{| \nabla u|}{\lambda} \right) \, dx < \frac{\varepsilon}{2\max \{ 1,\lambda^{n} \}}
\end{equation}
if  $\delta\leq\delta'$. Combining inequalities \eqref{feb29}, \eqref{feb27} and \eqref{feb32} yields
\begin{equation}\label{feb33}
 \sum_i \psi\big(d(u(Q_i))\big)   < c\delta \max \{ 1,\lambda^{n} \}+\frac{\ep}{2} 
\end{equation}
provided $\delta \leq \delta'$. Therefore,  if  $\delta < \min\big\{\delta',\frac{\ep}{2c(1+\lambda^{n})}\big\}$, then,
owing to inequality   \eqref{jan26}, one has that
$$
\H_{\omega(\sigma)}^{\psi}\bigl( u(E)  \bigr) 
\leq \sum_{i} \psi \bigl( d(u(Q_{i})) \bigr) < \ep
$$
for all $\sigma \in (0,1)$. Hence,   property \eqref{feb6} follows, since $\lim_{\sigma \to 0^+}\omega (\sigma )=0$.
\\
It remains to show that this property continues to hold for an arbitrary set $E \subset \Omega$. For any set $E$ of this kind, 
define the increasing sequence of sets $\{E_j\}$ as 
\begin{equation}\label{Ek}
E_j = \bigl\{ x \in E : \, |x| < j \mbox{ and } \mathrm{dist}(x,\partial \Omega ) > \tfrac{1}{j} \bigr\}
\end{equation}
 for $j \in \mathbb{N}$.
Clearly each set $E_j$ has a compact  closure  in  $\Omega$,  and $E = \bigcup_{j} E_j$.
Given $\varepsilon >0$, let $\delta >0$ be such that property \eqref{feb6} holds for every set whose  closure is compact and contained in $\Omega$. Assume that $\H^{\vp}(E)< \delta$.  Hence, $\H^{\vp}(E_j)< \delta$ as well, and, by property \eqref{feb6} applied with $E$ replaced by $E_j$, we have that $\H^{\psi}\big( u(E_{j})\big)<\varepsilon$.  Inasmuch as  $u(E) = \bigcup_{j} u(E_{j})$, from property \eqref{subad}
we can conclude that $\H^{\psi}\big( u(E)\big)\leq  \lim_{j\to \infty} \H^{\psi}\big( u(E_{j})\big)\leq \varepsilon$.  The proof of   \eqref{feb6} is complete.

\smallskip
\par\noindent Let us next consider Part (ii). Assume that $\H^{\vp}_{\infty}(E) <\delta$, where the number $\delta >0$ will be fixed later.
Hence, there exists $\sigma >0$ such that
$\H^{\vp}_{\sigma}(E)<\delta$ as well.  Furthermore, owing to property \eqref{apr1},  we have that $\Lambda_{2\sigma}^{\vp}(E) \leq c_{n}\H^{\vp}_{\sigma}(E)$. Thus, by  equation 
\eqref{feb26}, there exists a constant $c=c(n,\kappa)$ such that  $\Lambda_{2\sigma}^{\vp_\kappa}(E)< c \delta$, and hence  we may select  a family $\{ Q_{i} \}$
of non-overlapping dyadic cubes such that $E \subset \bigcup_{i}Q_i$, 
\begin{equation}\label{Feb70}
\sup _i d(Q_{i}) < 2\sigma,
\end{equation}
and
\begin{equation}\label{Feb50}
\sum_{i} \vp \bigl(\kappa d(Q_{i}) \bigr) < c\delta.
\end{equation}
Inequalities \eqref{feb28} and \eqref{feb29} continue to hold also in this case, whence  
\begin{align}\label{Feb73}
 \sum_i \psi\big(d(u(Q_i))\big) 
\leq \max \{ 1,\lambda^{n} \}\bigg( \sum_i \varphi (\kappa d(Q_i)) + \kappa^{n+1}   \int_{\bigcup_i Q_i} A\bigg(\frac{|\nabla u|}\lambda\bigg)\, dx\bigg),
\end{align}
provided that $\lambda > \|\nabla u\|_{L^A(\Rn)}$.
Also, the same chain as in equation \eqref{feb31} tells us that
\begin{align}\label{Feb72}
  \Le^n \Big(\bigcup_i Q_i\Big) \leq \frac {c\delta}{n^{\frac{n}{2}}}
  \frac{ \sigma ^{n}}{\vp \bigl( \kappa  \sigma \bigr)} \leq \frac {c\delta}{n^{\frac{n}{2}}\vp \bigl( \kappa\bigr)}.
\end{align}
Consequently, there exists $\delta ' = \delta' (n,m,A, \varphi, u,  \varepsilon)>0$ such that
\begin{equation}\label{Feb74}
\kappa^{n+1}\int_{\bigcup_{i}Q_i} \! A\left( \frac{| \nabla u|}{\lambda} \right) \, dx < \frac{\varepsilon}{2\max \{ 1,\lambda^{n} \}}
\end{equation}
if $\delta < \delta '$. On making use   inequalities \eqref{Feb73} and \eqref{Feb74}  one can infer that
$$\H_{\infty}^{\psi}\bigl( u(E)  \bigr) 
\leq \sum_{i} \psi \bigl( d(u(Q_{i})) \bigr) < \ep.
$$
Property \eqref{Feb60} is thus established.
\qed

\medskip
\par\noindent
\textit{Proof of Corollary \ref{cor}.} 
Let $E$ be any subset of $\Omega$ such that $\H^\vp(E)=0$.
 Consider any increasing sequence  $\{\Omega_j\}$ 
 of open subsets of $\Omega$, such that $\overline {\Omega_j} \subset \Omega$ for $j \in \N$, and  $\Omega = \bigcup_{j} \Omega_j$. 
 Define, for each  $j \in \N$, the set $E_j = E \cap \Omega_j$. Obviously, $\H^\vp(E_j)=0$ for every  $j \in \N$. 
 Since $u \in {\rm V}^{1,A}_{\rm loc}(\Omega, \Rm)$, one has that 
$u \in {\rm V}^{1,A}(\Omega_j, \Rm)$ for every $j \in \N$. An application of Theorem \ref{main}, Part (i), with $\Omega$ and $E$ replaced by $\Omega_j$ and $E_j$, respectively, 
tells us that $\H^{\psi}\bigl( u(E_{j}) \bigr) = 0$ for every $j \in \N$. Inasmuch as $\{u(E_j)\}$ is an increasing sequence of sets, and
 $u(E) = \bigcup_{j} u(E_{j})$, we deduce through property \eqref{subad} that $\H^{\psi}( u(E)) = 0$.
\qed
 
\medskip
\par\noindent
\textit{Proof of Theorem \ref{thm:measdist}}
The main step consists in proving that, if $u \in {\rm V}^{1,A}(\Omega, \Rm)$, then 
\begin{equation}\label{measdistort1}
\H^{\psi}\bigl( u(E) \bigr) \leq c_{n}\lim_{r \to 0^{+}}   \Theta_{\psi}^0\bigl( \Theta_{B^{-1}}^\infty(r)\bigr)   \Theta_{\psi}^0\bigl( \kappa \| \nabla u \|_{\LL^{A}( \Omega )}^{+} \bigr) \H^{\vp}(E),
\end{equation}
for every set $E \subset \Omega$. This will evidently establish Part~(ii).
\\ In order to prove inequality \eqref{measdistort1},
 we may clearly assume that  $\H^{\vp}(E)< \infty$. We also assume in addition, for the time being, that $\overline E$ is compact   and  $\overline E \subset \Omega$, and we choose   $\varrho_0$  such that the set   $B_{\varrho_0}(E)$, defined as in \eqref{Brho},  
is compactly contained in $\Omega$. Our assumptions on $A$ ensure that the function $u$ is uniformly continuous in $B_{\varrho_0}(E)$, with a modulus of continuity $\omega$, say.
\\
Now, fix $ \varepsilon  >0$. In view of \eqref{feb25},   for each $\delta \in (0, \varrho_0 )$ we can select a family $\{ Q_{i} \}$ of non-overlapping dyadic cubes such that
$E \subset \bigcup_{i}Q_{i}$, $d \bigl( Q_{i} \bigr)< \delta$ and
\begin{equation}\label{apr30}
\sum_{i} \vp \bigl( d(Q_{i})\bigr) < c_n \H^{\vp}\bigl( E \bigr) + \varepsilon .
\end{equation}
We may   also  assume, without loss of generality, that $Q_i \subset B_{\varrho_0}(E)$.    Next, fix any    $r \in (0,1)$.  Thanks to inequality \eqref{modifyfeb16}, one  can deduce that 
\begin{equation}\label{Feb40}
\vp \left( J_{r}^{-1}\left( \frac{d\bigl( u(Q_{i}) \bigr)}{\kappa \lambda}\right) \right) \leq \vp \bigl( d(Q_{i}) \bigr) + 
\frac{d(Q_{i})^n}{r}B\left( \frac{d\bigl( u(Q_{i})\bigr)}{\kappa \lambda d(Q_{i})}\right)
\end{equation}
for every  $i$.
Since $d\bigl( u(Q_{i}) \bigr) \leq \omega( \delta )$,    by inequality \eqref{FebA} we may choose $\delta < \varrho_0$ such
that, if $d\bigl( u(Q_{i}) \bigr) > 0$, then
\begin{equation}\label{Feb41}
\vp \left( J_{r}^{-1}\left( \frac{d \bigl( u(Q_{i}) \bigr)}{\kappa \lambda} \right) \right) >
\frac{\psi \left( \frac{d( u(Q_{i}))}{\kappa \lambda}\right)}{\Theta_{\psi}^0\bigl( \Theta_{B^{-1}}^\infty(r)^{+} \bigr) + \varepsilon}.
\end{equation}
On the other hand, all cubes for which $d\bigl( u(Q_{i}) \bigr) = 0$ can be disregarded in estimating $ \H^{\psi}(u(E))$, since $\psi(0)=0$. On decreasing $\delta$, if necessary, we infer from inequality \eqref{Feb41} and the definition of the function $\Theta_{\psi}^0$ that
\begin{equation}\label{Feb42}
\vp \left( J_{r}^{-1}\left( \frac{d\bigl( u(Q_{i}) \bigr)}{\kappa \lambda}\right) \right) > 
\frac{\psi \bigl( d\bigl( u(Q_{i}) \bigr) \bigr)}{\Theta_{\psi}^0\bigl( \Theta_{B^{-1}}^\infty(r)^{+} \bigr)+\varepsilon} \, \frac{1}{\Theta_{\psi}^0( \kappa \lambda ) +\varepsilon }
\end{equation}
for every index $i$ such that $d\bigl( u(Q_{i}) \bigr) > 0$.
Combining inequalities \eqref{Feb40} and \eqref{Feb42} with   inequality \eqref{poincare} yields 
\begin{equation}\label{Feb56}
\psi \bigl( d\bigl( u(Q_{i}) \bigr) \bigr) \leq \big(\Theta_{\psi}^0( \kappa \lambda )+ \varepsilon\big)\biggl( \Theta_{\psi}^0\bigl( \Theta_{B^{-1}}^\infty(r)^{+} \bigr) + 
\varepsilon \biggr)\left( \vp \bigl( d(Q_{i}) \bigr) +\frac{\kappa}{r}\int_{Q_i} \! A \left( \frac{| \nabla u|}{\lambda} \right) \, \dd x \right)
\end{equation}
for every index $i$ such that $d\bigl( u(Q_{i}) \bigr) > 0$.
Summing this inequality over all these indices $i$ and recalling inequality \eqref{apr30} tell us that
\begin{equation}\label{Feb47}
\H^{\psi}_{\omega ( \delta )}\bigl( u(E) \bigr) \leq \big(\Theta_{\psi}^0( \kappa \lambda )+ \varepsilon\big)\biggl(
\Theta_{\psi}^0 \bigl( \Theta_{B^{-1}}^\infty(r)^{+} \bigr) + \varepsilon \biggr)\left( c_{n}\H^{\vp}(E) +\varepsilon +
\frac{\kappa}{r}\int_{\bigcup Q_i} \! A \left( \frac{| \nabla u|}{\lambda} \right) \, \dd x \right) .
\end{equation}
We claim that the integral in this inequality tends to $0$ as $\delta \to 0^+$.
Indeed,
\begin{eqnarray*}
         \Le^{n}\Big( \bigcup_{i} Q_{i}\Big) &=& n^{-\frac{n}{2}}\sum_{i}d(Q_{i})^n =
n^{-\frac{n}{2}}\sum_{i}\vp \bigl( d(Q_{i}) \bigr) 
\frac{d(Q_{i})^{n}}{\vp \bigl( d(Q_{i}) \bigr)}\\
&
\leq & n^{-\frac{n}{2}}\sum_{i}\vp \bigl( d(Q_{i} ) \bigr)
\frac{\delta^{n}}{\vp ( \delta )} \leq n^{-\frac{n}{2}}\bigl( c_{n}\H^{\vp}(E)+\varepsilon \bigr)\frac{\delta^{n}}{\vp ( \delta )},
\end{eqnarray*}
where the first inequality holds thanks to property \eqref{feb8}. Our claim hence follows, owing to assumption \eqref{Feb4}.
Consequently,  from inequality \eqref{Feb47} we obtain that
$$
\H^{\psi}\bigl( u(E) \bigr) \leq \big(\Theta_{\psi}^0( \kappa \lambda )+ \varepsilon\big)\biggl(
\Theta_{\psi}^0 \bigl( \Theta_{B^{-1}}^\infty(r)^{+} \bigr) +\varepsilon \biggr) \bigl(
c_{n}\H^{\vp}(E)+\varepsilon \bigr).
$$
Inequality \eqref{measdistort1} follows, owing to the arbitrariness of $r\in (0, 1]$, of
 $\varepsilon > 0$, and of $\lambda > \| \nabla u \|_{\LL^{A}( \Omega )}$.
\\ Our next task is to show that inequality  \eqref{measdistort1} continues to hold for an arbitrary subset of $\Omega$. Of course, we can still assume that $\H^\vp(E)<\infty. $
Consider the increasing sequence of sets $\{E_j\}$ defined in \eqref{Ek}.  
Note that $\overline {E_j}$ is a compact set contained in  $\Omega$,  and $E = \bigcup_{j} E_j$. Inequality  \eqref{measdistort1} can thus be applied with $E$ replaced by $E_j$ for each $j \in \N$, and tells us that
\begin{align}\label{Feb48}
\H^{\psi}\bigl( u(E_{j}) \bigr) & \leq  c_{n}\Theta_{\psi}^0 \bigl( \kappa \| \nabla u \|_{\LL^{A}( \Omega )}^{+}
\bigr) \lim_{r\to 0^+}\Theta_{\psi}^0\bigl( \Theta_{B^{-1}}^\infty(r)^{+} \bigr)\H^{\vp}(E_{j}) \\ \nonumber & \leq  c_{n}\Theta_{\psi}^0 \bigl( \kappa \| \nabla u \|_{\LL^{A}( \Omega )}^{+}
\bigr)  \lim_{r\to 0^+}\Theta_{\psi}^0\bigl( \Theta_{B^{-1}}^\infty(r)^{+} \bigr) 
\H^{\vp}(E).
\end{align}
Inasmuch as  $u(E) = \bigcup_{j} u(E_{j})$, inequality \eqref{measdistort1} for the set $E$ follows thanks to property \eqref{subad}.
Part~(ii) is thus fully proved.
\\ Part (i) is a straightforward consequence of inequality  \eqref{measdistort1} under the assumption that $u \in {\rm V}^{1,A}(\Omega, \Rm)$. It remains to show that the latter assumption can be relaxed by just requiring that $u \in {\rm V}^{1,A}_{\rm loc}(\Omega, \Rm)$. This can be accomplished again by an invasion argument for $E$ by subsets.  Specifically, 
suppose that $E$ is any subset of $\Omega$ such that $\H^\vp(E)<\infty$. Let $\{ \Omega_j\}$ be an increasing sequence 
 of open subsets of $\Omega$, whose closure is also contained in $\Omega$ and such that $\Omega = \bigcup_{j} \Omega_j$. Set $E_j = E \cap \Omega_j$ for $j \in \N$ and note that $E_j$ is a subset of $\Omega_j$ with  $\H^\vp(E_j)<\infty$.  Furthermore, $u \in {\rm V}^{1,A}(\Omega_j, \Rm)$. Owing to inequality \eqref{measdistort1}, applied with $\Omega$ and $E$ replaced by $\Omega_j$ and $E_j$, and to assumption \eqref{scalingassump}, one has that
$\H^{\psi}\bigl( u(E_{j}) \bigr) = 0$ for every $j \in \N$. Assertion \eqref{Feb1} hence follows via property \eqref{subad}, since
 $u(E) = \bigcup_{j} u(E_{j})$.
\qed

\noindent
\textit{Proof of Theorem \ref{thm:contdist}.}
We may assume that $\H^{\vp}_{\infty}(E)>0$, the case when 
$\H^{\vp}_{\infty}(E)=0$ being
covered by Theorem \ref{thm:measdist}. Let us suppose, for the time being, 
 that $\H^{\vp}_{\infty}(E) < \infty$.
Thanks to  the second inequality in \eqref{apr1}, one has that $\Lambda_{\infty}^{\vp}(E) \leq c_{n}\H^{\vp}_{\infty}(E)$.
Thus, given any $\varepsilon > 0$, there exists a family $\{ Q_{i} \}$
of non-overlapping dyadic cubes such that $E \subset \bigcup_{i}Q_i$ and
\begin{equation}\label{Feb64}
\sum_{i} \vp \bigl( d(Q_{i}) \bigr) < c_{n}\H^{\vp}_{\infty}(E) + \varepsilon .
\end{equation}
From inequalities \eqref{Feb40}, \eqref{FebB}, the definition of the function $\Theta_{\psi}$
and 
inequality \eqref{poincare} we deduce (in analogy with inequality \eqref{Feb56}) that, for each index $i$,
\begin{equation}\label{Feb63}
\psi \bigl( d\bigl( u(Q_{i}) \bigr) \bigr) \leq \big(\Theta_{\psi}(\kappa \lambda){ \color{blue} + \varepsilon}\big) \bigl( \Theta_{\psi}\bigl( \Theta_{B^{-1}}(r)^{+}\bigr) + \varepsilon \bigr)\left( \vp \bigl( d(Q_{i}) \bigr) +
\frac{\kappa}{r}\int_{Q_i} \! A\left( \frac{| \nabla u|}{\lambda} \right) \, \dd x \right)
\end{equation}
for every   $r>0$ and
$\lambda > \| \nabla u \|_{\LL^{A}( \Rn )}$.
Summing over all indices $i$ in inequality \eqref{Feb63}, and making use of \eqref{Feb64} enable one to infer that
\begin{equation*}
\H^{\psi}_{\infty} \bigl( u(E) \bigr) 
\leq  \big(\Theta_{\psi}(\kappa \lambda)+ \varepsilon\big)\bigl( \Theta_{\psi}\bigl( \Theta_{B^{-1}}(r)^{+}\bigr) + \varepsilon \bigr)\left( c_{n}\H^{\vp}_{\infty}(E) + \varepsilon +
\frac{\kappa}{r}\int_{\Rn} \! A\left( \frac{| \nabla u|}{\lambda} \right) \, \dd x \right).
\end{equation*}
Now we let $\lambda \to \| \nabla u \|_{\LL^A(\Omega)}^{+}$ and $\varepsilon \to 0^+$ to deduce that
$$
\H^{\psi}_{\infty} \bigl( d\bigl( u(E) \bigr) \bigr) \leq \Theta_{\psi}\bigl( \kappa \| \nabla u \|_{\LL^A}^{+} \bigr)\Theta_{\psi}
\left( \Theta_{B^{-1}}(r)^{+}\right) \biggl( c_{n}\H^{\vp}_{\infty}(E)+\frac \kappa r \biggr).
$$
Choosing $r = \tfrac{\kappa}{c_n\H^{\vp}_{\infty}(E)}$ in this inequality yields inequality \eqref{Feb54}.
\\  It remains to consider the case when $\H^{\vp}_{\infty}(E)= \infty$. One can then just apply inequality \eqref{Feb54} to each set $E_j$ of  the sequence   $\{E_j\}$ defined by $E_j = E \cap \{x\in \rn: |x|<j\}$ for $j \in \N$. Passing to the limit as $j \to \infty$ and making use of property \eqref{subcont} yields the conclusion also in this case. \qed
 
 \noindent
 \textit{Proof of Theorem \ref{nlusin}.} 
 As a first step, we establish inequality \eqref{Feb60-t37} for 
 sets $E$ whose closure is compact and contained in $\Omega$. 
 Fix any such set $E$ and set  $\delta = \mathrm{dist}\bigl( E, \partial \Omega \bigr)$.  
 Let $\lambda > \| \nabla u \|_{L^{A}( \Omega )}$. Given any $\varepsilon \in (0,\delta )$, we can select a family $\bigl\{ Q_i \bigr\}$ of non-overlapping dyadic cubes such that 
 $E \subset \bigcup_{i} Q_i \subset B_{\varepsilon}(E) \subset   \overline{B_{\varepsilon}(E)}\subset \Omega$, $d(Q_{i}) < \varepsilon$ and $\Le^{n}\bigl( \bigcup_{i}Q_i \bigr) < \Le^{n}(E)+\varepsilon$.  
 Here, $B_{\varepsilon}(E)$ denotes the set defined as in \eqref{Brho}. Let $r>0$ be a number to be specified later. The function $J_r$, associated with $B$ and with the function 
 $\varphi (s)= s^n$ as in \eqref{Jr}, obeys
\begin{equation}\label{march1}
J_r(s) = s B^{-1}(r) \quad \text{for $s \geq 0$.}
\end{equation}
Thus, from inequalities \eqref{modifyfeb16} and \eqref{poincare} one infers that
 \begin{align*} 
 \frac{r}{B^{-1}(r)^{n}}\bigg( \frac{d \bigl( u(Q_{i}) \bigr)}{\kappa \lambda}\bigg)^{n} \leq r d(Q_{i})^{n}+d(Q_{i})^{n}B\bigg( \frac{d\bigl( u(Q_{i})\bigr)}{\kappa \lambda d(Q_{i})}\bigg)\leq  r d(Q_{i})^{n}+\kappa \int_{Q_i} \! A \left( \frac{| \nabla u|}{\lambda}\right) \, d x  
 \end{align*}
for every $i$. Hence,
 \begin{align}\label{march2}
d \bigl( u(Q_{i}) \bigr)^n  &\leq \kappa^{n}\lambda^{n}B^{-1}(r )^{n}\bigg( \Le^{n}(Q_{i})n^{\frac{n}{2}}+\frac \kappa r \int_{Q_i} \! A \left( \frac{| \nabla u|}{\lambda}\right) \, d x\bigg)
 \end{align}
for every $i$.
 Summing this inequality over all indices $i$ yields
 \begin{align*}
 \Le^{n}\bigl( u(E) \bigr) &\leq  \kappa^{n}\lambda^{n}B^{-1}(r)^{n}\bigg( (\Le^{n}(E)+\varepsilon )n^{\frac{n}{2}}+ \frac{\kappa}{r}\int_{\bigcup_{i}Q_i} \! A \left( \frac{| \nabla u|}{\lambda} \right) \, dx \bigg)\\
 &\leq \kappa^{n}\lambda^{n}B^{-1}(r)^{n}\left( (\Le^{n}(E)+\varepsilon )n^{\frac{n}{2}}+ \frac{\kappa}{r} \right),
 \end{align*}
 whence, by the arbitrariness of $\lambda > \| \nabla u \|_{L^{A}( \Omega )}$ and of $\varepsilon \in (0, \delta )$, one deduces that
 \begin{equation}\label{march3}
 \Le^{n}\bigl( u(E) \bigr) \leq \bigl( \kappa \| \nabla u \|_{L^{A}( \Omega )} B^{-1}( r ) \bigr)^{n} \left( \Le^{n}(E)n^{\frac{n}{2}}+\frac{\kappa}{r} \right).
 \end{equation}
 Inequality \eqref{Feb60-t37} follows from \eqref{march3}, by either choosing
$r= \frac{\kappa}{n^{\frac{n}{2}}\Le^{n}(E)}$ or letting $r\to 0^+$, according to whether $\Le^{n}(E)>0$ or $\Le^{n}(E)=0$. In the latter case, the fact that $\lim_{r\to 0^+} \Phi_B(r)=0$ plays a role.
\\ Assume now that $E$ is an arbitrary subset of $\Omega$ and let $E_j$ be the set defined, for $j\in \N$, as in \eqref{Ek}.
The sequence of sets $\{E_j\}$ is increasing and consists of sets with compact closure contained in $\Omega$. Furthermore, $E = \bigcup E_j$. 
We already know that inequality  \eqref{Feb60-t37} holds  with $E$ replaced by each $E_j$. Namely,
 $$
\H^{n}_{\infty}\bigl( u(E_j) \bigr) \leq 2\kappa^{n+1} \| \nabla u\|_{\LL^{A}( \Omega )}^{n} \Phi_{B}\big( \Le^{n}(E_j)n^{\frac{n}{2}}/\kappa\big)
%
%
 $$
for every $j \in \N$.
Since the function $\Phi_{B}$ is increasing, this inequality implies that
$$
\H^{n}_{\infty}\bigl( u(E_j) \bigr) \leq 2\kappa^{n+1} \| \nabla u\|_{\LL^{A}( \Omega )}^{n} \Phi_{B}\big( \Le^{n}(E)n^{\frac{n}{2}}/\kappa\big)
%
%
 $$
 for every $j\in \N$. Inasmuch as  $\bigcup_{j} u(E_{j}) = u(E)$,   property \eqref{subcont}
ensures that $\lim _{j \to \infty}\H^{n}_{\infty}\bigl( u(E_{j}) \bigr)=\H^{n}_{\infty}\bigl( u(E) \bigr)$.  Inequality \eqref{Feb60-t37} is fully proved.
\qed

 \section{Sharpness}\label{S:cex}

This final section is devoted to discussing the 
 sharpness of the conclusions of Theorems~\ref{main}  for some examples  from Section~\ref{S:main}.

\paragraph*{I.}  Let $A$ and $\varphi$ be as in  Example~\ref{ex4}. Namely, let $A$ be any Young function such that $i_\infty (A)>n$ and let $\varphi$ be a gauge function fulfilling
$\liminf_{r \to 0^+} \frac{ \varphi (r^\gamma)}{ \varphi (r)} > 0$   for every $\gamma >0$.
We have shown that
$$\psi (r) \sim \varphi (r) \quad \text{near zero.}$$ 
Thus the result is sharp, since this is the best possible choice of the gauge function $\psi$ one can hope for.   In particular, this includes the case exhibited in equation  \eqref{sh1}, i.e. maps $u\in \VV^{1,p}(\mathbb R^n, \mathbb R^n)$\, with $p>n$,
and gauge functions satisfying $\varphi (r) \sim \big(\log \tfrac 1r\big)^\beta$  near zero, for some   $\beta <0$.

\paragraph*{II.}
As displayed in  equation \eqref{kr0}, when  the classical Sobolev space $ \VV^{1,p}(\mathbb R^n, \mathbb R^n)$, with $p>n$, and the standard Hausdorff measure $\H^\alpha$, with
$\alpha \in (0,n]$,  are in question, a~special instance of Example~\ref{ex1} reproduces Kaufman's theorem. In Kaufman's paper, the exponent in the Hausdorff
measure $\H^{\frac{\alpha p}{p+\alpha-n}}$ is shown to be optimal. 
\\
The next result augments this conclusion and ensures that the gauge function $\psi (r) = r^{\frac{\alpha p}{p+\alpha-n}}$, which defines the Hausdorff measure $\H^{\frac{\alpha p}{p+\alpha-n}}$ in the target space, does not even
admit an improvement in the scale of powers of a~logarithm, at least if these powers are not too small.

 \begin{theorem}\label{main-ex-s2}  {\sl Assume that $n\in \mathbb N$,  $p>n$ and $\alpha \in (0,n)$. Set $\sigma= \frac{\alpha p}{p+\alpha-n}$ and let $\gamma >1+\sigma$. 
Let $\phi$ be a gauge function such that
\begin{equation}\label{jan30}
\phi(r) \sim r^\sigma\big(\log \tfrac 1r\big)^\gamma  \quad \text{near zero.}
\end{equation}
Then there exist a  mapping $u\in \VV^{1,p}(\mathbb R^n, \mathbb R^n)$ and a compact set $E\subset\R^n$ such that
\begin{equation}\label{ex-f-cs2}
0< \H^{\alpha }(E) < \infty \qquad \mbox{\rm and }\qquad \H^{\phi}(u(E))=\infty.
\end{equation}}
\end{theorem}
 
\paragraph*{III.}
Here we focus on maps from Example \ref{ex2}, dealing with Young functions
$$\text{ $A(t)$  equivalent to $t^n(\log t)^{q} $ near infinity,}$$
in the special case when the gauge function 
$$ \varphi (r) \sim \big(\log\tfrac1r\big)^{\beta} \quad \text{ near zero,}$$
 for some  $\beta <0$.
The result presented in that example asserts that the corresponding gauge function  in the target space is of merely power type, and fulfills
$$ \psi (r) \sim r^{\frac {n\beta}{n-q-1+\beta}} \quad \text{ near zero.}$$
The following theorem tells us that not only is this gauge function the best possible among all standard Hausdorff measures of power type, but it is also
optimal in the scale of \lq\lq powers times powers of a logarithm" type gauge functions,   at least if the power of the logarithm is not too small.  

\begin{theorem}\label{main-ex-s}  {\sl Assume that $n\in \mathbb N$, 
\begin{equation}\label{jan32}
A(t)\sim  t^n(\log t)^{q} \quad \text{\rm near infinity } 
\end{equation}
for some  $q>n-1$, and 
\begin{equation}\label{ex-k4}
\varphi(r) \sim\big(\log\tfrac1r\big)^{\beta} \quad \text{near zero}
\end{equation}
for some $\beta <0$. Set $\sigma=\frac {n\beta}{n-q-1+\beta}$, let  $\mu >1+\sigma$, and let  
$\phi$ be a gauge function such that
 \begin{equation}\label{jan33}
\phi(r) \sim r^\sigma\big(\log \tfrac 1r\big)^\mu \quad \text{near zero.}
\end{equation}
Then there exist a mapping $u\in \VV^{1,A}(\Rn , \Rn )$ and a compact set $E\subset\Rn$ such that
\begin{equation}\label{ex-f-cs1}
0< \H^{\varphi }(E) < \infty \qquad \mbox{\rm and }\qquad \H^{\phi}(u(E))=\infty.
\end{equation}}
\end{theorem}


The proofs of Theorems \ref{main-ex-s2} and \ref{main-ex-s} are reminiscent of arguments from \cite{Kau} and their modification in \cite[section~4.2]{hencl-h2015}.
In particular, they involve a construction resting upon probabilistic techniques. Additional  difficulties arise in the present situation,
which are specific of the Orlicz ambient space. Theorem \ref{main-ex-s} will be established in detail. The proof of Theorem \ref{main-ex-s2}
follows along the same lines, and it is in fact simpler, since the Orlicz-Sobolev space agrees with a classical Sobolev space in this case.
For brevity, this proof will be omitted. 

The question of whether the lower bounds on the exponents $\gamma$ and $\mu$ in Theorems \ref{main-ex-s2} and \ref{main-ex-s}, respectively, can be removed is open. This is a natural  interesting problem, whose solution seems however to require approaches and techniques substantially different from those employed in the present paper and in previous contributions in the literature. An even more challenging question is of course whether $\psi$ is best possible among all gauge functions, for any given $A$ and $\varphi$.

\medskip

In preparation for the proof of Theorem  \ref{main-ex-s}, 
we recall a few definitions and establish  a few preliminary lemmas.
\\
The Hausdorff dimension $\dim_\H(E)$ of a set $E \subset \rn$ is defined as 
$$\dim_\H(E)
=\inf\{\alpha>0:\H^\alpha(E)=0\}.$$
One also has that $\dim_\H(E)
=\sup\{\alpha>0:\H^\alpha(E)=\infty\}$.
\\ Given a  Radon measure $\mu$ on a~compact set $E\subset\R^n$ and a number $\alpha >0$, we set 
$$I_\alpha(\mu) =\int\limits_E\int\limits_E|x-y|^{-\alpha}\,d\mu(x)\,d\mu(y).$$
It is well known that 
\begin{equation}\label{cde-1}
\dim_\H(E)=\sup\{\alpha>0: \text{there exists $\mu\ne0$ such that $I_\alpha(\mu)<\infty$}\},
\end{equation}
see, e.g., \cite[Theorem~8.9]{Mat}.
\\ The following result (see \cite[Theorem~8.7]{Mat}) will also be needed in our proof of Theorem  \ref{main-ex-s}.

\begin{lemma}\label{cglem}  Let $\phi$ be a gauge function, and let $E$ is compact set in $\rn$ . Assume that there exists a nonnegative Radon measure $\mu$ such that $0<\mu(E)<\infty$ and  
$$
\int\limits_E\int\limits_E\frac{\dd \mu(x)\,\dd \mu(y)}{\phi \bigl( |x-y| \bigr)} <\infty.
$$
Then $\H^\phi(E)=\infty$.
\end{lemma}

The following probabilistic lemma extends, with analogous proof,   \cite[Lemma~4.4]{bal2013}. 
In its statement, $\mathbf E$ denotes expectation, and $\|\{a_i\}\|_\infty$  the $\ell^\infty$ norm of a sequence $\{a_i\}$. Henceforth, $\log$ stands for $\log_2$.

\begin{lemma}\label{prob-lem} Let $\{X_i\}_{i=1}^\infty$ be a countable sequence of independent random variables,
identically distributed according to the uniform distribution on the unit ball~$B$ in $\R^n$.  Let $\phi : [0, \infty) \to [0, \infty)$ be the function defined by
$$\phi (r) = r^\sigma\log^\mu\bigl(b+\frac1r\bigr) \quad \text{for $r \geq 0$,}$$ 
with  $\sigma \in (0, n)$, $\mu\in\R$, and $b$ large enough for the function $\phi$ to be increasing.  Assume that 
 $\{a_i\}\in\ell_1$.  Then there exists a constant $c$, depending on $n$, $\sigma$, $\mu$, and $b$, such that 
 \begin{equation}\label{prob-est}
  \mathbf E\biggl[\biggl(\phi\biggl(\biggr|\sum\limits_{i=1}^\infty a_i X_i\biggr|\biggr)\biggr)^{-1}\biggr]\le c\, \frac1{\phi\bigl(\|\{a_i\}\|_\infty\bigr)}
\end{equation}
\end{lemma}

Let us denote by $L^n\log^q L(\rn)$ the Orlicz space built upon the Young function 
$A(t)= t^n\log^q(t+2)$, with $q>n-1$, and equipped with the Luxemburg norm. The corresponding Orlicz-Sobolev space will be denoted by $\VV ^1L^n\log^2 L(\rn)$.
\\ We use the notation $Q$ for a cube in $\rn$  with sides parallel to the coordinate axes, and $\ell(Q)$ for its sidelength. By 
$2Q$ we denote  the cube with double sidelength and the same center as $Q$. 
Finally, $Q(r)$ stands for  the cube, centered at the origin,  with $\ell(Q)=r$. 

\begin{lemma}\label{scal-lem}  For every $j\in\N$, there exists a function $\eta_{j}\in \VV ^1L^n\log^2 L(\rn)$ such that 
$$0 \leq \eta (x) \leq 1 \quad \text{for $x \in \rn$,}$$
$$\eta_j(x)=\begin{cases}1 \qquad \text{if \, $x\in Q(2^{-2^{j}})$}
\\[10pt]
0\qquad \text{if \, $x\notin Q(2^{-2^{j-1}})$,}
\end{cases}
and 
$$
\begin{equation}\label{ex-6}
  \|\nabla \eta_j\|_{L^n\log^q L(\rn)}\le c\,2^{j(q-n+1)/n},
\end{equation}
where $c$ is a constant independent of~$j$.
\end{lemma}
\textit{Proof.}
Let us define, for $j\in \N$,
$$
\eta_j(x)=\begin{cases}1\qquad & \text{if $2|x|_\infty\le2^{-2^{j}}$} 
\\
\frac1{2^{(j-1)}}\biggl(-2^{ (j-1)}+\log\frac1{2|x|_\infty}\biggr) \qquad  &\text{if $2^{-2^{j}}<2|x|_\infty<2^{-2^{(j-1)}}$}
\\
0 \qquad & \text{if $2|x|_\infty\ge 2^{-2^{(j-1)}}$,}
\end{cases}
$$
where $|x|_\infty=\max\limits_{i=1,\dots,n}|x_i|$ and   $x=(x_1,x_2,\dots,x_n)\in\R^n$. 
In particular, given $r>0$,
$$x\in Q(r) \quad \text{if and only if} \quad 2|x|_\infty\le r.$$ 
Thus, the function $\eta_j$ satisfies the first two properties claimed in the statement.
\\ As for inequality \eqref{ex-6}, notice that
$$
|\nabla\eta_j(x)| \begin{cases} \sim \frac1{2^j\,|x|} \quad & \text{if $2^{-2^{j}}<2|x|_\infty<2^{-2^{(j-1)}}$ }
\\ 
=0  \quad & \text{otherwise.}
\end{cases}
$$
 Now, set
 $\lambda =2^{j(q-n+1)/n}$. Thus, on passing to polar coordinates, one  obtains that
\begin{eqnarray*}
\int\limits_{\R^n}A\bigg(\frac{|\nabla\eta_j|}{\lambda}\bigg)\,dx &\leq& \frac{c}{2^{jn}\lambda^n}\int\limits_{2^{-2^{j}}}^{2^{-2^{(j-1)}}}\frac1r\log^q\Big(2+\frac1{\lambda2^jr}\Big)\, \dd r\\
&\le& \frac{c}{2^{j(q+1)}}\int\limits_{2^{-2^{j}}}^{2^{-2^{(j-1)}}}\frac1r\log^q\Big(\frac1r\Big)\, \dd r \le \frac{c}{2^{j(q+1)}}2^{j(q+1)}\le c.
\end{eqnarray*}
for some constant $c$ independent of $j$. Hence, inequality \eqref{ex-6} follows.
 \hfill \qed

\begin{lemma}\label{scal-lem2}{\sl  Let $N \in \N$ and let $v_{i} \colon \Rn \to\R$, $i=1,\dots,N$, be a  family of  functions in $\VV L^n\log^q L(\rn)$, with disjoint supports,  and such that
$$
\|\nabla v_1\|_{L^n\log^q L(\rn)}=\|\nabla v_2\|_{L^n\log^q L(\rn)}=\dots=\|\nabla v_N\|_{L^n\log^q L(\rn)}.
$$
Then 
\begin{equation}\label{apr40}
\biggl\|\sum\limits_{i=1}^N\nabla v_i\biggr\|_{L^n\log^q L(\rn)}\le N^\frac1n\|\nabla v_1\|_{L^n\log^q L(\rn)}.
\end{equation}
}
\end{lemma}

\noindent
\textit{Proof.}
Without loss of generality, we can assume that $\|\nabla v_i\|_{L^n\log^q L(\rn)}=1$ for $i=1,\dots,N$, namely,
$$
\int\limits_{\Rn} \! |\nabla v_i|^n\log^q\bigl(|\nabla v_i|+2\bigr)\, \dd x=1\qquad \text{for $i=1,\dots,N$.}
$$
On setting $\lambda =N^\frac1n$, one has that
\begin{eqnarray*}
  \int\limits_{\Rn} \! A\biggl(\frac1\lambda\biggr|\sum\limits_{i=1}^N\nabla v_i\biggr|\biggr) \, \dd x &=&
  \frac1{\lambda ^n}\sum\limits_{i=1}^N\int\limits_{\R^n}|\nabla v_i|^n\log^q\biggl(\frac1\lambda|\nabla v_i|+2\biggr)\,\dd x\\
                    &\le&
\frac1{\lambda^n}\sum\limits_{i=1}^N\int\limits_{\R^n}|\nabla v_i|^n\log^q\bigl(|\nabla v_i|+2\bigr)\, \dd x=\frac{N}{\lambda^n}=1.
\end{eqnarray*}
Hence, inequality  \eqref{apr40} follows.
\qed

\medskip
\par\noindent
We are now ready to provide a proof of Theorem \ref{main-ex-s}.

\medskip
\par\noindent
\textit{Proof of Theorem \ref{main-ex-s}}.
For simplicity of notation, throughout the proof we rename the parameter $\beta$, and set 
$\nu =-\beta$. Hence, $\nu >0$ and  
\begin{equation}\label{exp-k4}
\varphi(r) =\bigl(\log\tfrac1r\bigr)^{-\nu} \quad \text{for $r>0$.}
\end{equation}
For each $j \in \N$, consider ~a system of dyadic cubes $\A_j$ in $[0,1]^n$  enjoying the following properties:  
$$\ell(Q)=2^{-2^j}\qquad  \text{for every  $Q\in\A_j$},$$ 
the cardinality $\#\bigl(\A_j\bigr)$ of $\A_j$ 
satisfies 
\begin{equation}\label{paramm-2}
\#\bigl(\A_j\bigr) = [2^{j\nu}],
\end{equation}
where  $[\,\cdot\,]$ denotes integer part, and 
$$\mbox{ for every }Q\in\A_{j+1}\quad\mbox{ there exists }\,\widetilde Q \in \A_j\quad\mbox{ such that }\,Q\subset\tilde Q.$$
The construction of such a family is elementary and standard. 
\\
Let us set  
$$\A=\bigcup\limits_{j\in\N}\A_j \quad \text{and} \quad M=\bigcap\limits_{j\in\N}\bigcup\limits_{Q\in\A_j}Q.$$
One can verify that
\begin{equation}\label{hM}
0<\H^\vp(M)<\infty.
\end{equation}
Moreover, we can assume that the cubes
are evenly distributed,  and hence that  there exists a~universal constant $c>0$ such that, for every $j \in \N$,
\begin{equation}\label{exp-k7}
  \text{ if $Q\in\A_j$, then $\H^\vp(Q\cap M)\le c\,2^{-j}$,}
\end{equation}
and
\begin{equation}\label{exp-k7-}
\text{if $\widetilde Q\in\A_j$, then  $\#\{Q\in\A_{j+1}:Q\subset \widetilde Q\}\le c$.}
\end{equation}
Set  \begin{equation}\label{exp-k5}
\sigma=\frac {n\nu}{q-n+1+\nu}, \end{equation}
and fix $\mu>1+\sigma$. Define the function $\phi :[0, \infty) \to [0, \infty)$ by
$$ \phi(r)= \begin{cases} r^\sigma\log^\mu\bigl(b+\tfrac1r)\quad &\text{if $r>0$}
\\ 0 \quad &\text{if $r=0$,}
\end{cases}$$
where the~constant $b$ is large enough for the~function $\phi$ to be increasing.  


For each $Q\in\A_j$ let $\eta_Q$ be a function defined as an appropriate
translate of the function $\eta_j$ from Lemma~\ref{scal-lem}, which fulfills the following properties:
\\
(i)   $\eta_Q$ is supported in the cube~$Q'\supset \supset Q$ with $\ell (Q') = 2^{-2^{j-1}}$; 
\\
(ii) $0 \le\eta_Q\le1$;
\\
(iii) $\eta_Q=  1$ on $Q$;
\\
(iv) $\eta_Q= 0$ outside $Q'$;
\\
(v) $\|\nabla \eta_Q\|_{L^n\log^q L(\rn)}\le c\,2^{j(q-n+1)/n}$.

\medskip
\par\noindent
Let $\xi = \{\xi_Q\}_{Q\in\A}$ be a countable sequence of points in the unit ball $B$ of $\Rn$. 
Fix  $\delta>1$ such that 
\begin{equation}\label{paramm-1}
 1+ \delta\sigma<\mu.
\end{equation}
For each $j\in\N$, define    the function $u_{\xi,j}: \rn \to \rn$ as
$$
u_{\xi,j}(x)=\frac1{j^\delta}2^{-\frac{j\nu}{\sigma}}\sum\limits_{Q\in\A_j}\eta_Q(x)\xi_Q \quad \text{for $x \in \rn$,}
$$
and the function $u_\xi : \rn \to \rn$ as
$$
u_\xi(x)=\sum\limits_ju_{\xi,j}(x) \quad \text{for $x \in \rn$.}
$$
We   claim that there exists a universal constant $c>0$  such that 
\begin{equation}\label{exp-k10}
 \#\bigl\{Q\in\A_{j}:(\supp \eta_Q)\cap(\supp \eta_{\widetilde Q})\ne\emptyset\bigr\}\le c
\end{equation}
for every $j\in\N$ and every  $\widetilde Q\in\A_j$.
Indeed, if $Q,\widetilde Q\in\A_j$ and $(\supp \eta_Q)\cap(\supp \eta_{\widetilde Q})\ne\emptyset$, then $Q$ and $\widetilde Q$ either lie
in the same cube or in adjacent cubes
from $\A_{j-1}$. Moreover, any dyadic cube has exactly $3^n-1$ adjacent dyadic cubes of the same size.
These pieces of information, combined with property~(\ref{exp-k7-}), enable one to complete the proof of~(\ref{exp-k10}).
\\
From  equations (\ref{paramm-2}), (\ref{exp-k10}) and   Lemmas~\ref{scal-lem}--\ref{scal-lem2} we deduce that there exists a constant $c$, independent of $\xi$ and $j$, such that
$$
\|\nabla u_{\xi,j}\|_{L^n\log^q L(\rn)}\le c\,\frac1{j^\delta}2^{-\frac{j\nu}{\sigma}}\ 2^{j\nu/n}\ 2^{j(q-n+1)/n}
=c\frac1{j^\delta}.
$$
Therefore,
$$
\|\nabla u_{\xi}\|_{L^n\log^q L(\rn)}\le c\sum\limits_j\frac1{j^\delta}<\infty.
$$
This shows that the mapping $u_\xi \in  \VV^1 L^n\log^q L(\rn)$, and hence, in particular,  
$u_\xi$ is continuous, since we are assuming that $q>n-1$.  Also, we have that $\supp u_\xi\subset [0,1]^n$. 
 
The remaining part of the proof is devoted to showing that,   for a generic choice
of~$\xi$,  the map~$u_\xi$ enjoys the property
\begin{equation}\label{exp-k-f}
\H^{\phi}\bigl( u_{\xi}(M) \bigr)=\infty.
\end{equation}
The desired conclusion will then follow by coupling equation \eqref{exp-k-f} with  \eqref{hM}.
\\ In order to establish property \eqref{exp-k-f}, let us regard $\xi=\{\xi\}_{Q\in\A}$
as a sequence of independent random variables,  which are identically distributed according to
the uniform probability distribution on  $B$.
 Consider
the measure $\mu_\xi=(u_\xi)_\#(\H^\vp \restrict M)$ supported on the compact set $u_\xi(M)$, i.e., $\mu_\xi$ is the push-forward of the measure $\H^\vp$ restricted to $M$ via the map
$u_\xi$. We claim that its expectation satisfies the condition
\begin{equation}\label{exp-k11}
  \mathbf E_\xi\biggl(\int\limits_{u_\xi(M)}\int\limits_{u_\xi(M)}\frac{\dd \mu_\xi(x)\,\dd \mu_\xi(y)}{\phi \bigl( |x-y| \bigr)} \biggr)<\infty.
\end{equation}
In order to prove inequality \eqref{exp-k11}, we make use Fubini's theorem and rewrite the left-hand side of inequality \eqref{exp-k11} as
\begin{equation}\label{apr44}
\int\limits_M\int\limits_M\mathbf E_\xi\biggl(\frac{1}{\phi \bigl(\bigl| u_\xi(x)-u_\xi(y)\bigr|\bigr)} \biggr)\,
\dd \H^\vp(x)\, \dd \H^\vp(y).
\end{equation}
We have that
$$
u_\xi(x)-u_\xi(y)=\sum\limits_{Q\in\A}a_Q(x,y)\xi_Q \quad \text{for $x, y \in \rn$,}
$$
where 
\begin{equation}\label{exp-k12}
  a_Q(x,y)=\frac1{j^\delta}2^{-\frac{\nu  j}\sigma}\bigl(\eta_Q(x)-\eta_Q(y)\bigr) \quad \text{ for $Q\in\A_j$.}
\end{equation}
The sequence  $\{a_Q\}$ clearly belongs to~$\ell^\infty$.
Thus, thanks to inequality  (\ref{prob-est}), inequality \eqref{apr44} will follow if we show that
\begin{equation}\label{apr45}
\int\limits_M\frac1{\phi\bigl(\|a(x,y)\|_\infty\bigr)}\, \dd \H^\vp(x)
\le c,
\end{equation}
for some constant~$c$ is independent of~$y$. Let us fix~$y\in M$. For each $x\in M$, let us denote
by~$j(x)$ the largest integer such that both~$x$ and~$y$ either lie in the same cube or in adjacent cubes
from~$\A_{j(x)}$. Hence, there exist cubes ~$\widetilde Q_1\ni x$ and $\widetilde Q_2\ni y$ such that $\widetilde Q_1,\widetilde Q_2\in\A_{j(x)+1}$ and $2\tilde Q_1\cap Q_2=\emptyset$. 
Moreover, there exist cubes 
$Q_1\ni x$ and $Q_2\ni y$ such that $Q_1,Q_2\in\A_{j(x)+2}$ and 
$$
\bigl(\supp\eta_{Q_1}\bigr)\cap\bigl(\supp\eta_{Q_2}\bigr)=\emptyset.
$$
Consequently, $\eta_{Q_1}(x) = 1$ and 
$\eta_{Q_1}(y) = 0$.  Therefore, there exists  a~nonzero coefficient $a_Q$
of size~$\frac{2^{-(j(x)+2)\frac\nu\sigma}}{(j(x)+2)^\delta}$.
Hence
$$
\|a(x,y)\|_\infty\ge\frac{2^{-j(x)\frac\nu\sigma}}{j(x)^\delta}.
$$
The monotonicity of the function $\phi$ guarantees that
\begin{equation}\label{paramm-4}
  \frac1{\phi\bigl(\|a(x,y)\|_\infty\bigr)}= \frac1{\|a(x,y)\|^\sigma_\infty\log^\mu\bigl(b+\frac1{\|a(x,y)\|_\infty}\bigr)} \le c\,j(x)^{\delta\sigma-\mu}\,\,2^{\nu j(x)}
\end{equation}
for some constant $c$.
\\
Now, for each 
$y\in M$ and   $j\in\N$,  there exist at most $3^n$ dyadic cubes $Q\in\A_{j}$ with the property that there exists  $x\in Q\cap M$
with $j(x)=j$. Thereby,  equation (\ref{exp-k7}) implies that
$$
\H^\vp\bigl(\bigl\{x\in M:j(x)=j\bigr\}\bigr)\le 3^n\H^\vp(Q\cap M)\le c \ 2^{-j\nu}.
$$
Finally, owing to inequality (\ref{paramm-4}) and assumption \eqref{paramm-1},
$$
\int\limits_M\frac1{\phi\bigl(\|a(x,y)\|_\infty\bigr)}\, \dd \H^\vp(x)\le
c\sum\limits_{j=1}^\infty2^{{\color{blue}-j\nu}}\biggl(j^{\delta\sigma-\mu}\,2^{\nu j}\biggr)=c\sum\limits_{j=1}^\infty j^{\delta\sigma-\mu}
<\infty.
$$
Inequality \eqref{apr45}, and hence property  \eqref{exp-k11}, are  thus established. 
\\ On the other hand, property  \eqref{exp-k11}
 clearly entails that, almost surely with respect to $\xi$,  
$$
\int\limits_{u_\xi(M)}\int\limits_{u_\xi(M)}\frac{\dd \mu_\xi(x)\,\dd \mu_\xi(y)}{\phi \bigl( |x-y| \bigr)}<\infty.
$$
Hence, equation \eqref{exp-k-f} follows  via Lemma~\ref{cglem}. The proof is complete.
\qed

%

 \section*{Compliance with Ethical Standards}\label{conflicts}

\smallskip
\par\noindent 
{\bf Funding}. This research was partly funded by:   
\\ (i) Research Project 201758MTR2  of the Italian Ministry of Education, University and
Research (MIUR) Prin 2017 ``Direct and inverse problems for partial differential equations: theoretical aspects and applications'';   
\\ (ii) GNAMPA   of the Italian INdAM - National Institute of High Mathematics (grant number not available).

\smallskip
\par\noindent
{\bf Conflict of Interest}. The authors declare that they have no conflict of interest.

\end{document}